\documentclass[10pt,twoside,english,a4paper]{article}
\usepackage{etex}
\usepackage{caption}
\usepackage[latin1]{inputenc}
\usepackage[intlimits] {amsmath}   % muss for defjs# stehen
\usepackage{graphicx}
\usepackage{psfrag}
\usepackage{ifthen}
\usepackage{fancyhdr}
\usepackage{rotating}
\usepackage{multirow}
\usepackage{booktabs}              % Dicke Tabellenlinien
\usepackage{amssymb}
\usepackage{amsbsy}
\usepackage{bm}                    % fette Schrift in Formeln
\usepackage{bbm}
\usepackage{babel}
\usepackage{theorem}
\usepackage{upgreek}  % gerade (roman) griechische Buchstaben z.B. \upalpha
\usepackage{pifont}   % zusaetzliche Schriften, z.B. eingekreiste Zahlen
\usepackage[active]{srcltx}   % aktive suche im xdv
\usepackage{subcaption}
\usepackage{suetterl}  % Suetterlin Schrift
\newfont{\suetdbl}{suet14 scaled 2000}  % Suetterlin skaliert um Faktor 2,000

\usepackage{yfonts}  % altdeutsche Schrifte und Initialen
\newfont{\gothdbl}{ygoth scaled 2000} 
\newfont{\frakdbl}{yfrak scaled 2000}
\newfont{\swabdbl}{yswab scaled 2000}
\usepackage{caption}
\usepackage{subcaption}
\usepackage{placeins} 
\usepackage[normalem]{ulem}  % zusaetzliche Unterstreichungsformen wie z.B.:
\usepackage[textsize=footnotesize,textwidth=1.7cm]{todonotes}
\usepackage[]{defjs3}      % eigene Definitionen: Option BoldVarsRoman
\usepackage{stmaryrd}
\usepackage{subcaption}
\usepackage[]{defjs3}  
\usepackage{amsfonts}
\usepackage{graphicx}
\usepackage{amsmath}
\usepackage{comment}
\usepackage{algorithm2e}
\RestyleAlgo{boxed}
\SetKwProg{Block}{}{}{}
\DontPrintSemicolon
\SetAlCapSkip{1em}
\usepackage[numbers,sort&compress]{natbib}
\bibpunct{[}{]}{,}{n}{}{;}
\fancypagestyle{plain}{%
  \fancyhead{}%
  \fancyfoot[c]{\sffamily\thepage}%
}
\makeatletter                           % Leere Fuellseiten
\def\cleardoublepage{\clearpage\if@twoside \ifodd\c@page\else
  \hbox{}
  \vspace*{\fill}
  \thispagestyle{empty}
  \newpage
  \if@twocolumn\hbox{}\newpage\fi\fi\fi}
\makeatother
\def\Curl{\operatorname{Curl}}      % Curl  ohne Argument
\def\curl{\operatorname{curl}}       % curl  ohne Argument
\newcommand{\Ce}{\mathbb{C}_{\mathrm{e}}}
\newcommand{\CeV}{\widetilde{\mathbb{C}}_{\mathrm{e}}}
\newcommand{\Cc}{\mathbb{C}_{\mathrm{c}}}
\newcommand{\Cmicro}{\mathbb{C}_{\mathrm{micro}}}
\newcommand{\Cmacro}{\mathbb{C}_{\mathrm{macro}}}
\newcommand{\CmicroV}{\widetilde{\mathbb{C}}_{\mathrm{micro}}}
\newcommand{\CmacroV}{\widetilde{\mathbb{C}}_{\mathrm{macro}}}

\newcommand{\Cmatrix}{\mathbb{C}_{\mathrm{matrix}}}

\newcommand{\Lc}{L_{\mathrm{c}}}
\newcommand{\IL}{\mathbb{L}}
\newcommand{\R}{\mathbb{R}}
\newcommand{\Bdis}{\bP}
\newcommand{\dis}{P}
\begin{document}
\unitlength1.0cm
\frenchspacing

%=== sections ===
\thispagestyle{empty}

\vspace{1mm}
\ce{\bf 
\Large A computational approach to identify the material }

\ce{\bf \Large parameters of the relaxed micromorphic model }

\vspace{3mm}
\ce{Mohammad Sarhil$^{1,2,\ast}$,  Lisa Scheunemann$^2$, Peter Lewintan$^{3,4}$,  J\"org Schr\"oder$^1$ and Patrizio Neff$^4$}

\vspace{3mm}
\ce{$^1$Institute of Mechanics, University of Duisburg-Essen}
\ce{Universit\"atsstr. 15, 45141 Essen, Germany}
\ce{\small e-mail: 
mohammad.sarhil@uni-due.de
}    

\vspace{2mm}
\ce{$^2$Chair of Applied Mechanics, RPTU Kaiserslautern}
\ce{Gottlieb-Daimler-Str., 67663 Kaiserslautern, Germany}

\vspace{2mm}
\ce{$^3$Faculty of Mathematics, Karlsruhe Institute of Technology}
\ce{Englerstr. 2, 76131 Karlsruhe, Germany}

\vspace{2mm}
\ce{$^4$Chair of Nonlinear Analysis and Modeling, Faculty of Mathematics, }
\ce{ University of Duisburg-Essen, Thea-Leymann-Str. 9,}
\ce{ 45141 Essen, Germany}

\vspace{2mm}
\begin{center}
{\bf \large Abstract}
\bigskip

{\small
\begin{minipage}{14.5cm}
\noindent
{We determine the material parameters in the relaxed micromorphic generalized continuum model for a given periodic microstructure in this work. This is achieved through a least squares fitting of the total energy of the relaxed micromorphic homogeneous continuum to the total energy of the fully-resolved heterogeneous microstructure, governed by classical linear elasticity. The relaxed micromorphic model is a generalized continuum that utilizes the $\Curl$ of a micro-distortion field instead of its full gradient as in the classical micromorphic theory, leading to several advantages and differences. The most crucial advantage is that it operates between two well-defined scales. These scales are determined by linear elasticity with microscopic and macroscopic elasticity tensors, which respectively bound the stiffness of the relaxed micromorphic continuum from above and below. While the macroscopic elasticity tensor is established a priori through standard periodic first-order homogenization, the microscopic elasticity tensor remains to be determined. Additionally, the characteristic length parameter, associated with curvature measurement, controls the transition between the micro- and macro-scales. Both the microscopic elasticity tensor and the characteristic length parameter are here determined using a computational approach based on the least squares fitting of energies. This process involves the consideration of an adequate number of quadratic deformation modes and different specimen sizes. We conduct a comparative analysis between the least square fitting results of the relaxed micromorphic model, the fitting of a skew-symmetric micro-distortion field (Cosserat-micropolar model), and the fitting of the classical micromorphic model with  two different formulations for the curvature; one simplified formulation involving only one single characteristic length and a simplified  isotropic curvature with three parameters. The relaxed micromorphic model demonstrates good agreement with the fully-resolved heterogeneous solution after optimizing only four parameters. The ``simplified" full micromorphic model, which includes isotropic curvature and involves the optimization of seven parameters, does not achieve superior results, while the Cosserat model exhibits the poorest fitting. 
 }  

\end{minipage}
}
\end{center}


\section{Introduction}

Architected materials or metamaterials are unconventional materials with exceptional mechanical properties that depend on the intricate geometry of the underlying complex microstructure, rather than the bulk properties of their constituent materials. They can be engineered to fulfill specific functionalities. Nonetheless, they often exhibit size-effects, meaning that their effective properties change with the variations in material size when scale separation does not hold.  Generally,  size-effects can involve both increasing and decreasing stiffness with reducing the size \cite{WheFraRic:2015:isa,KirAmsHue:2023:otq}.  In this work, we consider only the case that smaller is comparatively stiffer. On the other hand, materials with intricate microstructures often require homogenization techniques since a complete resolution of the underlying microstructure is typically infeasible for usual engineering problems due to the associated substantial computational costs. First-order homogenization methods \cite{Zoh:2004:hma,ZohWri:2005:ait} and standard elasticity theories, by their nature, are incapable of capturing size-effects phenomena. Generalized continua are enhanced models which are capable of incorporating the size-effects without the need to account for the fully-resolved microstructures. Two families of generalized continua are commonly found in the literature. The first family extends the kinematics to include additional degrees of freedom such as an independent ``micro-" rotation, stretch, strain or full deformation field, such as those found in the classical Eringen-Mindlin  micromorphic theory \cite{Eri:1968:mom,SuhEri:1964:nto,Eri:1973:tom} and the Cosserat theory \cite{CosCos:1909:tof,NefJeoMueRam:2010:lce}.  The second family of generalized continua incorporates higher-order differential operators in the energy or motion functional, such as strain- or stress-gradient \cite{ForBer:2020:sge}, as seen in gradient elasticity models \cite{Aif:2011:otg,AskAif:2011:gei}.   A classification of generalized continua constructed by a micromechanical approach is available in \cite{AlaGanRedSad:2023:hog}.

\begin{figure}[htpb]
\center
	\unitlength=1mm
	\begin{picture}(140,65)
	\put(20,2){\def\svgwidth{12 cm}{\small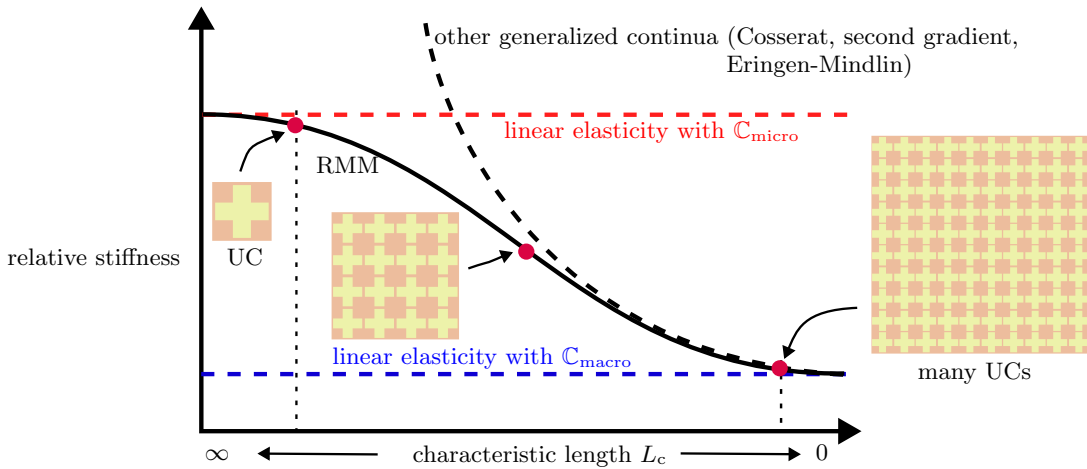}}
	\end{picture}
	\caption{The stiffness of the relaxed micromorphic model (RMM) is bounded from above and below. Other generalized continua exhibit unbounded  stiffness for small sizes. For large values of the characteristic length, linear elasticity with a micro elasticity tensor is recovered (one UC) while linear elasticity with a macro elasticity tensor is obtained for small values of  the characteristic length (many UCs).   }
	\label{Figure:intro}
\end{figure}

The relaxed micromorphic model (RMM) is a generalized continuum description which simplifies the form of the assumed strain energy compared to the classical micromorphic theory by using a relaxed curvature in terms of the $\Curl$ of a micro-distortion field rather than its full gradient \cite{NefGhiMadPlaRos:2014:aup,MadNefGhiRos:2016:rat,DemRizColNefMad:2022:uem}. Utilizing only the $\Curl$ of the micro-distortion field offers several advantages, notably reducing the number of material parameters. The well-posedness for the important case of symmetric force stress is proven using new generalized Korn's inequality  \cite{Lewintan2021K,Lewintan:2021:tfg,GemLewNef:2023:cov,GmeLewNef:KMS:2024}.  An important characteristic  of the relaxed micromorphic continuum is its bounded stiffness from below and above, allowing the material parameters to be related to two well-defined scales \cite{GouRizLewBerSkyMadNef:2023:gff} which is impossible for the classical micromorphic or gradient elasticity theories. This feature establishes the relaxed micromorphic model as a linear elasticity model operating between these two scales: the microscopic scale, described by linear elasticity with a micro elasticity tensor representing the maximum stiffness exhibited by the assumed metamaterial, and the macroscopic scale, characterized by linear elasticity with a macro elasticity tensor obtained using standard periodic first-order homogenization methods where scale separation holds. In the RMM, the characteristic length parameter plays a critical role in scaling correctly with the specimen size and controlling the transition between the micro- and macro-scales. Therein, large values of the characteristic length correspond to zooming into the ``stiff" microstructure, for example, a unit-cell (UC), while small values result in an effective ``soft" classical homogeneous response for large structures. Figure \ref{Figure:intro} illustrates how this unique behavior distinguishes the relaxed micromorphic model from other generalized continua which exhibit unbounded stiffness for arbitrarily small specimens (i.e., large values of the characteristic length). In previous attempts, the micro elasticity tensor was defined as the stiffest response at the unit-cell level. In \cite{NefEidMad:2019:ios}, focusing on band gaps, the micro elasticity tensor was determined by the L\"owner matrix supremum of elasticity tensors under affine Dirichlet conditions. However, it was found to be too soft for size-effects in the bending regime \cite{SarSchNefSch:2022:mts}. This prompted an expansion of our understanding of the micro elasticity tensor, particularly by incorporating non-affine Dirichlet conditions \cite{SarSchSchNef:2023:oti, SarSchSchNef:2023:seo}, resulting in a micro elasticity tensor calibrated specifically for bending. However, a homogenization procedure for the identification of all the unknown parameters, including the characteristic length, has not yet been established.

The identification of material parameters for enriched continua is a highly non-trivial task which remains largely unresolved despite many attempts in the literature. Various methods have been proposed for the homogenization of fully-resolved heterogeneous microstructures into the Cosserat continuum in \cite{ForSab:1998:com, Hue:2019:otm, RedAlaNasGan:2021:htc}, different variants of the gradient elasticity continuum in  \cite{SarYidAba:2023:coh,AbaBar:2021:ami, AbaYanPap:2019:aca, BacPagDalBig:2018:iof, KhaNii:2020:asg, LahGodGan:2022:sis, SchKruKeiHes:2022:cho, SkrEre:2020:ote, Wee:2021:nho, YanMue:2021:seo, YanAbaTimMue:2020:dom, YanAbaMueetal:2022:voa}, and the classical Eringen-Mindlin micromorphic continuum in  \cite{AlaGanRedSad:2021:com, BisPoh:2017:amc, For:2002:hma, Hue:2022:iom, Hue:2017:hoa, RokAmePeeGee:2019:mch, RokAmePeeGee:2020:emc, RokZemDosKry:2020:ris, ZhiPohTayTan:2022:dfm}. However, these methods have yet to yield a universally accepted solution. Many fundamental questions arise in the context of homogenization towards generalized continua, most of which have been addressed within the framework of the standard first-order homogenization theory \cite{TriJanAufDieFor:2012:eof,GanWazRed:2023:fih}. The definition of a representative volume element and the choice of the boundary condition are the first obvious issues which one faces in the homogenization towards higher-order continua. While in first order homogenization the condition of continuity of the local fields at the interface
of unit-cell results in a periodicity requirement of the micro-displacement fluctuation
field, see for example \cite{VonHelSchSch:2024:otr}, this periodicity requirement becomes more or less irrelevant in the absence of scale separation and an overall strain gradient loading in the framework of higher-order homogenization, \cite{BacGam:2010:soh}. Consequently, the higher order moduli are dependent on the choice of the representative volume element. Alternatively, the analysis can be done on a cluster of unit-cells to get rid of edge effects and considering the converged behavior in the central unit-cell but some zero-energy modes are obtained \cite{ForTri:2011:gca}. Another crucial point is that the average second gradient cannot be strictly controlled by the boundary condition, and usually, multiple modes are triggered simultaneously. Consequently, the selection of the relevant higher-order polynomial coefficients becomes very complex. This can be fixed by volumetric  constraints \cite{RalSte:2012:mlc}. In this study, we aim to circumvent the numerous previously unanswered questions by employing a non-classical homogenization strategy to determine the remaining unknown parameters of the relaxed micromorphic model. This strategy is based on the least squares fitting of the energy of the homogeneous relaxed micromorphic continuum with that of fully discretized specimens and does not require the use of classical (or non-classical) micro-macro transition relations. By considering various deformation modes, whether random or not, and different specimen sizes, we identify the unknown parameters of the relaxed micromorphic model. A pertinent methodology is outlined in \cite{AbaYanPap:2019:aca}, which selects $n$ deformation modes to determine $n$ unknowns and solves the resulting $n$ equations precisely.

The paper follows the outline: In Section \ref{sec:model}, we revisit the relaxed micromorphic model, exploring its strong and weak form, alongside the related boundary conditions and the modified curvature which is scaled with the number of considered unit-cells. In Section \ref{section:unit_cell}, we detail the metamaterial geometry and the material parameters of the unit cell. Section \ref{section:linearelasticity} is dedicated to presenting an algorithm that serves as a motivational example and a conceptual validation. This algorithm defines the stiffness matrix of an equivalent homogeneous continuum for a single unit-cell under both affine Dirichlet and periodic boundary conditions. The algorithm is then further extended in Section \ref{section:RMM} to encompass the case of an equivalent relaxed micromorphic homogeneous continuum, accompanied by numerical examples. Additionally, we compare the fitting results of the relaxed micromorphic model, the Cosserat model and the classical micromorphic model with a simplified curvature and an  isotropic curvature.  We draw our conclusions in Section \ref{Section:Conclusions}. For this work,  we limit our consideration to the planar case, in which the isotropic curvature energy in terms of $\Curl \Bdis$ has only one free parameter \cite{DagMatLewBerDanNef:2024:}. We confine the optimization to align with the assumed cubic unit cell, which is not a limitation, as we precisely know the anisotropy properties of the relaxed micromorphic model.

\sect{\hspace{-5mm}. The relaxed micromorphic model}
\label{sec:model}  
The relaxed micromorphic model (RMM) is a generalized continuum. Each material point's kinematics are described, similar to the general micromorphic theory \cite{EriSub:1964:nto,Min:1964:msi,SuhEri:1964:nto}, by a standard displacement vector $\bu\colon\B\subseteq\R^3\to\R^3$ and a non-symmetric micro-distortion field  $\Bdis\colon\B\subseteq\R^3\to\R^{3\times3}$. The displacement and the micro-distortion fields are defined by minimizing the energy functional  
\begin{equation}
\label{eq:pot}
\Pi (\bu,\Bdis)= \int_\B W\left(\nabla \bu,\Bdis,\Curl \Bdis \right) \, - { \overline\bbf}\cdot{\bu} \, \,\textrm{d}V\ - \int_{\partial \B_t}  \overline\bt \cdot \bu  \,\textrm{d}A  \longrightarrow\ \min\,,
\end{equation}
with $(\bu,\Bdis)\in H^1(\B)\times H(\Curl,\B)$. The vector $ \overline\bbf$  describes the applied body force. The vector $ \overline\bt$ is the traction vector acting on the boundary $\partial \B_t \subset \partial \B$. The elastic energy density $W$ reads 
\begin{equation}
\label{eq:W}
\begin{aligned}
W\left(\nabla \bu,\Bdis,\Curl \Bdis \right) =  \frac{1}{2} \Big( & \symb{ \nabla \bu - \Bdis} : \Ce : \symb{ \nabla \bu - \Bdis}  +   \sym \Bdis : \Cmicro: \sym \Bdis  \\
& + \skewb{ \nabla \bu - \Bdis} : \Cc : \skewb{ \nabla \bu - \Bdis}  +  \mu \, \Lc^2 \, \textrm{Curl} \Bdis : \IL :\textrm{Curl} \Bdis  \Big)\,.
\end{aligned}
\end{equation}
Here, $\Cmicro,\Ce > \bzero $ are fourth-order positive definite standard elasticity tensors, $\Cc \ge \bzero$ is a fourth-order positive semi-definite rotational coupling tensor, $\IL$ is a positive definite fourth-order tensor acting on non-symmetric arguments.  The  characteristic length parameter ($\Lc > 0 $) is related to the size of the microstructure. It allows to scale the number of considered unit-cells keeping all remaining parameters of the model scale-independent. The macro-scale with $\Cmacro$ and the micro-scale with $\Cmicro$ are retrieved for the limiting cases $\Lc \rightarrow 0$ and  $\Lc \rightarrow \infty$, respectively, if suitable boundary conditions are applied, see \cite{NefEidMad:2019:ios,SchSarSchNef:2022:lhb}.  The macro-scale elasticity tensor $\Cmacro$ associated with $\Lc \rightarrow 0$ is 
a priori uniquely defined by standard first-order periodic homogenization (the scale separation holds) while the micro-scale elasticity tensor $\Cmicro$ associated with $\Lc \rightarrow \infty$ represents the stiffest extrapolated response (zooming in the microstructure). The shear modulus $\mu$ in the curvature term has been added for dimensional consistency and is not necessarily related to $\Cmicro$ or $\Cmacro$.  The constitutive coefficients are assumed constant with the following symmetries 
\begin{equation}
\begin{aligned}
&(\Cmicro)_{ijkl} = (\Cmicro)_{klij} = (\Cmicro)_{jikl} \,, \qquad  &(\Cc)_{ijkl} = (\Cc)_{klij} \,, \\
&(\Ce)_{ijkl} = (\Ce)_{klij} = (\Ce)_{jikl} \,, \qquad   &{ (\IL)_{ijkl} = (\IL)_{klij}  \,,}
\end{aligned}
\end{equation} 
where $\Cmicro$ and $\Ce$ are connected to $\Cmacro$ through a Reuss-like homogenization relation { \cite{BabMadDagAbrGhiNeff:2017:taf}} (tilde stands for Voigt notation) with

\begin{equation}
\label{eq:reuss_like}
\CmacroV^{-1} = \CmicroV^{-1} + \CeV^{-1} \quad \Leftrightarrow \quad \CeV = \CmicroV \cdot  (\CmicroV-\CmacroV)^{-1} \cdot  \CmacroV\,.
\end{equation}

In order to build a link between the heterogeneous fully-detailed metamaterial  and the homogeneous relaxed micromorphic model, we address the size-effect by incorporating the number of considered unit-cells within the computational domain of  the reference heterogeneous material.  Thus, the modified energy functional reads

\begin{equation}
\begin{aligned}
W\left(\nabla \bu,\Bdis,\Curl \Bdis \right) =  \frac{1}{2} ( & \symb{ \nabla \bu - \Bdis} : \Ce : \symb{ \nabla \bu - \Bdis}  +   \sym \Bdis : \Cmicro: \sym \Bdis  \\
& + \skewb{ \nabla \bu - \Bdis} : \Cc : \skewb{ \nabla \bu - \Bdis}  +  \mu \, (\frac{\Lc}{n})^2 \, \textrm{Curl} \Bdis : \IL :\textrm{Curl} \Bdis  )\,.
\end{aligned}
\end{equation}

\begin{figure}[htpb]
\center
	\unitlength=1mm
	\begin{picture}(150,60)
	\put(0,0){\def\svgwidth{14 cm}{\small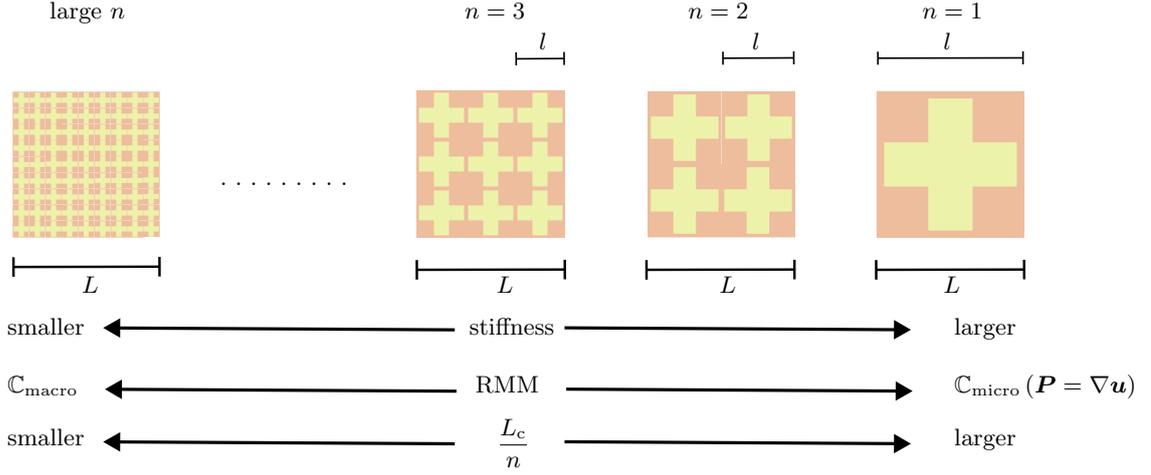}}
	\end{picture}
\caption{A depiction of the scaling $\Lc/n$. This scaling delivers the intended ``smaller is stiffer'' effect for  computations on a domain of fixed size and constant $\Lc$.}
\label{Figure:illustration:Lc/n}
\end{figure}

The scalar $n$ has been introduced resembling the number of the unit-cells considered for the heterogeneous metamaterial in the computational domain $\B = [-\dfrac{L}{2},\dfrac{L}{2}] \times  [-\dfrac{L}{2},\dfrac{L}{2}]  =  [-\dfrac{n \, l}{2},\dfrac{n \, l}{2}] \times  [-\dfrac{n \, l}{2},\dfrac{n \, l}{2}] $ as depicted in Figure \ref{Figure:illustration:Lc/n}.  The characteristic length parameter $\Lc$ is set then constant for any considered size.  This scaling of curvature by $\Lc/n$  is not ad hoc, but follows from a rigorous scaling argument,  cf. \cite{NefEidMad:2019:ios} which leads to the intended stiffening (smaller is stiffer). These considerations are not exclusive to the relaxed micromorphic model but  apply to the classical micromorphic, Cosserat and gradient elasticity models as well.

The variation of the potential with respect to the displacement yields the weak form
\begin{equation}
\begin{split}
\delta_{\bu} \Pi =& \int_{\B}\{  \underbrace{ \Ce : \symb{ \nabla \bu - \Bdis}+  \Cc : \skewb{ \nabla \bu - \Bdis}}_{\textstyle=:\Bsigma} \} :  {\nabla \delta \bu} -  \overline\bbf \cdot \delta \bu  \, \textrm{d}V - \int_{\partial \B_t}  \overline\bt \cdot \delta \bu  \, \textrm{d}A = 0\,,
\end{split}
\end{equation} 
which leads, using  integration by parts and employing the divergence theorem, to
 \begin{equation}
 \label{Eq:balance_linear}
\delta_{\bu} \Pi =  \int_{\B}   \{ \div \Bsigma +  \overline\bbf \} \cdot \delta \bu \, \textrm{d}V  = 0 \, ,
\end{equation}  
where  $\Bsigma$ is the non-symmetric force stress tensor. However, a symmetric force stress, i.e. $\Cc  = \bzero$, is always permitted, e.g when the consistent boundary condition is enforced on the boundary \cite{SarSchSchNef:2023:seo,SarSchNefSch:2022:mts}. In a similar way, the variation of the potential with respect to the micro-distortion field $\Bdis$ leads to the weak form
\begin{equation}
 \delta_{\Bdis} \Pi =   \int_{\B} \{ \Bsigma  - \underbrace{\Cmicro : \sym \Bdis}_{\textstyle =:\Bsigma_\textrm{micro} } \} :  \delta {\Bdis}  - \{ \underbrace{\mu \, \left(\dfrac{\Lc}{n}\right)^2  \IL : \Curl \Bdis}_{\textstyle =:\bbm} \} : \Curl \delta {\Bdis} \, \textrm{d} V  = 0\,, 
\end{equation}
which can be rewritten, using integration by parts and applying Stokes' theorem, as 
\begin{equation}
 \label{Eq:balance_angular}
 \delta_{\Bdis} \Pi =   \int_{\B}  \{ \Bsigma - \Bsigma_\textrm{micro}  - \Curl \bbm \} :  \delta {\Bdis} \, \textrm{d} V +  \int_{ \partial \B}  \{ \sum_{i=1}^3  \left( \bbm^i  \times \delta {\Bdis}^i  \right) \cdot \bn  \}  \;\; \textrm{d}A  
 = 0\,, 
\end{equation}
where the stress measures $\Bsigma_\textrm{micro}  $ and $\bbm$ are the micro- and moment stresses, respectively, $\bn$ is the outward unit normal vector on the boundary, and $\bbm^i$ and $ \delta {\Bdis}^i$ are the row vectors of the related second-order tensors. 
The strong forms of the relaxed micromorphic model, obtained from Equations  \ref{Eq:balance_linear} and \ref{Eq:balance_angular}, represent the generalized balance of linear and angular momentum  which  read with  the  associated boundary conditions 
\begin{equation}
 \begin{aligned}
  \div \Bsigma + \overline\bbf = \bzero  \quad &\textrm{on} \quad  \B ,  \qquad \qquad  \qquad \Bsigma - \Bsigma_\textrm{micro}  - \Curl \bbm  = 0  \quad  &\textrm{on}& \quad  \B \,, \\   
  \bu = \overline{\bu} \quad &\textrm{on} \quad  \partial \B_u\,, \qquad  \qquad \qquad \sum_{i=1}^3   \Bdis^i \times \bn =  \overline\bt_{P}  \quad &\textrm{on}& \quad  \partial \B_{P} \,, \\  
  \overline\bt = \Bsigma \cdot \bn \quad &\textrm{on} \quad  \partial \B_t \,, \qquad \qquad \qquad \sum_{i=1}^3   \bbm \bm^i \times \bn = \bzero \quad &\textrm{on}& \quad  \partial \B_m\,, 
   \end{aligned}
\end{equation}
 where $\partial \B_\dis \cap \partial \B_m = \partial \B_u \cap \partial \B_t  = \emptyset $ and  $\partial  \B_\dis \cup \partial \B_m = \partial \B_u \cup \partial \B_t = \partial \B $ (more details in \cite{SchSarSchNef:2022:lhb}).  By substituting the generalized balance of angular momentum into the generalized balance of linear momentum, a resulting, but not independent, balance equation reads
 
 \begin{equation}
  \div \Bsigma_\textrm{micro}  + \overline\bbf = \bzero \quad  \textrm{on} \quad  \B \,,
 \end{equation}

which does not appear in the Eringen-Mindlin micromorphic theory or the Cosserat (micropolar) model.  In the RMM, unlike the classical theory, only the tangential projection of the micro-distortion field can be described on the boundary, rather than the entire field. One obvious option is to link the tangential components of the micro-distortion $\Bdis$ and the deformation gradient $\nabla \bu$. This boundary condition, called {\bf consistent coupling condition}, was proposed in \cite{NefEidMad:2019:ios} and subsequently considered in \cite{SkyNeuMueSchNef:2021:CM,RizHueMadNef:2021:aso3,RizHueMadNef:2021:aso4,DagRizKhaLewMadNef:2021:tcc}. It reads 
 \begin{equation}
  \Bdis \cdot \Btau = \nabla \bu \cdot \Btau \, \Leftrightarrow \, {\Bdis}^i  \times \bn  =   {\nabla \bu}^i  \times \bn \quad \textrm{for} \quad i=1,2,3  \quad \textrm{on}\quad  \partial \B_\dis = \partial \B_u  \,,
 \end{equation}

  where $\Btau$ is the tangential vector on the boundary and ${\Bdis^i}$ and $\nabla {\bu}^i$ are the row-vectors of the associated tensors.  We may extend this boundary condition to parts of ${\partial} \B_m$ by enforcing the consistent coupling condition on $\partial \B_{\widehat{m}} \subseteq \partial \B_{{m}}$ using a penalty approach, as follows:
 \begin{equation}
 \Pi \Leftarrow \Pi + \int_{ \partial \B_{\widehat{m}}}   \frac{\kappa}{2} \sum_{i=1}^3 ||(\Bdis^i - \nabla {\bu}^i) \times \bn )||^2 \, \textrm{d}A  \,,
 \end{equation}
 where $\kappa$ is the penalty parameter.  The micro-distortion field for the three-dimensional case has the following general form
 \begin{equation}
\Bdis = \left[ \begin{array}{c}
(\Bdis^{1})^T \\
(\Bdis^{2})^T \\
(\Bdis^{3})^T \\
\end{array}\right] = \left[ \begin{array}{c c c}
 \dis_{11}  &  \dis_{12}  & \dis_{13}  \\   
 \dis_{21}  &  \dis_{22}  & \dis_{23}  \\ 
 \dis_{31}  &  \dis_{32}  & \dis_{33}  \\ 
\end{array}\right] \quad \textrm{with} \quad \Bdis^{i} = \left[ \begin{array}{c}
\dis_{i1} \\
\dis_{i2} \\ 
\dis_{i3}
\end{array} \right]\,  \textrm{for } i=1,2,3\,.
 \end{equation} 
We let the Curl operator act on the row vectors of the micro-distortion field $\Bdis$ as
 \begin{equation}
\Curl \Bdis = \left[ \begin{array}{c}
(\curl \Bdis^{1})^T \\
(\curl \Bdis^{2})^T \\
(\curl \Bdis^{3})^T \\
\end{array}\right] = \left[ \begin{array}{c|c|c}
 \dis_{13,2} - \dis_{12,3}	& \dis_{11,3} - \dis_{13,1}  &  \dis_{12,1} - \dis_{11,2} \\
 \dis_{23,2} - \dis_{22,3}	& \dis_{21,3} - \dis_{23,1}  &  \dis_{22,1} - \dis_{21,2} \\
 \dis_{33,2} - \dis_{32,3}	& \dis_{31,3} - \dis_{33,1}  &  \dis_{32,1} - \dis_{31,2} 
\end{array}\right] \,.
 \end{equation}

Various finite element formulations of the relaxed micromorphic model were presented for different cases: plane strain in \cite{SchSarSchNef:2022:lhb, SarSchNefSch:2021:oat}, antiplane shear in \cite{SkyNeuMueSchNef:2021:CM}, and the three-dimensional case in \cite{SkyNeuMueSchNef:2022:pam, SkyMueRizNeff:2023:hob}. A conforming finite element formulation for an even further relaxed curvature can be found in \cite{SkyNeuLewZilNEf:2024:nhc}. For the two-dimensional scenario, the micro-distortion field only has four non-vanishing components within the plane, and its Curl operator reduces to just two components out of the plane, namely $(\Curl \Bdis)_{13}$ and $(\Curl \Bdis)_{23}$,  such that 

 \begin{equation}
 \label{Eq:Curl_2D}
\Bdis = \left[ \begin{array}{c}
(\Bdis^{1})^T \\
(\Bdis^{2})^T \\
\bzero^T 
\end{array}\right] = \left[\begin{array}{c c c }
 \dis_{11}  &  \dis_{12}   & 0\\   
 \dis_{21}  &  \dis_{22}   & 0\\
 0 & 0 & 0 
\end{array}\right] \quad \textrm{and} \quad 
\textrm{Curl} \, \Bdis = 
\left[\begin{array}{c|c|c}
0	&   0  &  \dis_{12,1} - \dis_{11,2} \\
0   &   0  &  \dis_{22,1} - \dis_{21,2} \\
0 & 0 & 0
\end{array}\right]\,.
 \end{equation}

For this work, we utilize an $H(\Curl)$-conforming finite element, denoted as Q2NQ2, which employs Lagrange-type shape functions of the second-order for the displacement field, denoted as Q2, and  a N\'ed\'elec formulation of first-kind and second-order \cite{Ned:1980:mfe,Ned:1986:anf} for the micro-distortion field. For more details regarding the derivation  of shape functions and the FEM-implementation, the reader is referred  to \cite{SchSarSchNef:2022:lhb}. The simulations presented in this paper are performed within AceGen and AceFEM programs. The interested reader is referred to \cite{KorWri:2016:aofem,Kor:2009:aof} for more details on the latter two.

\section{The unit-cell and the material parameters} 
\label{section:unit_cell}
\begin{figure}[htpb]
    \centering
    \begin{minipage}{.3\textwidth}
        \centering
       	\begin{picture}(3,4)
      \put(0,0){\includegraphics[width=4cm]{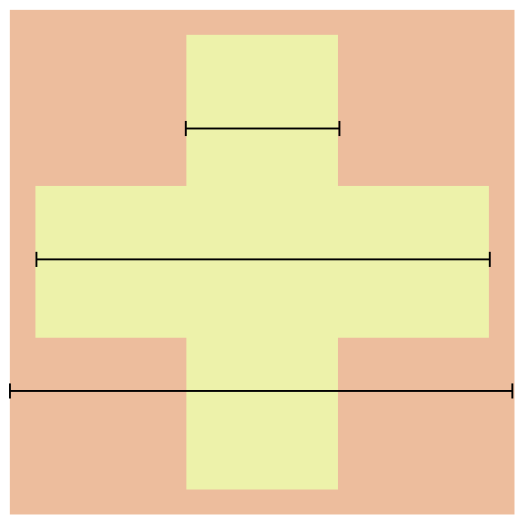}} 
      \put(1.838053286723413,1.20023909705){\color[rgb]{0,0,0}\makebox(0,0)[lb]{\smash{$l$}}}%
      \put(1.838053286723413,2.20023909705){\color[rgb]{0,0,0}\makebox(0,0)[lb]{\smash{$l_1$}}}%
       \put(1.838053286723413,3.20023909705){\color[rgb]{0,0,0}\makebox(0,0)[lb]{\smash{$l_2$}}}%
        \end{picture}      
    \end{minipage}
    \hfill
    \begin{minipage}{.65\textwidth}
        \centering
        \begin{tabular}{ccc}
            \toprule
           Lam\'e parameters & $\lambda$ [kN/mm$^2$] & $\mu$ [kN/mm$^2$]\\
            \midrule
            Matrix & $51.08 \,  $ & $26.32 \,  $  \\
            Inclusion & $51.08 \, 10^{-4} \, $ & $26.32  \, 10^{-4} \,  $  \\
            \bottomrule
        \end{tabular}
        
        \vspace{0.5 cm}
        \begin{tabular}{ccc}
            \toprule
            $l$ [mm] & $l_1$ [mm] & $l_2$ [mm]\\
            \midrule
           $1/n$  & $0.9 \, l$ & $0.3 \, l$ \\
            \bottomrule
        \end{tabular}
    \end{minipage}
            \caption{Unit-cell with the material and geometrical parameters.}
            \label{Fig;unit_cell}
\end{figure}
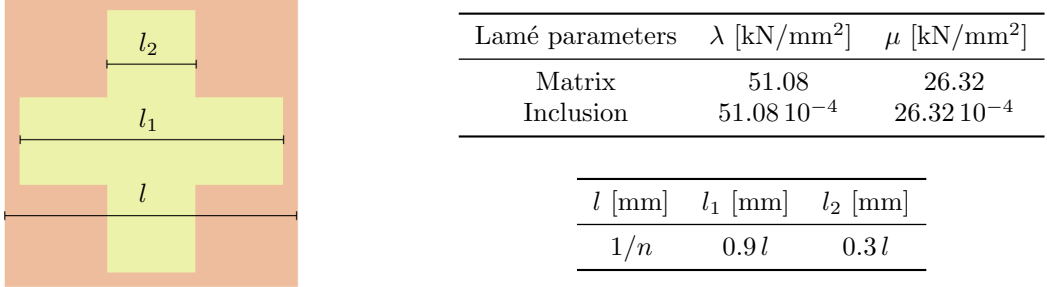

In the following, we consider  a unit-cell consisting of a stiff matrix (aluminum) and a swiss-cross shape ultra-soft inclusion where the Lam\'e parameters differ by a factor $10000$. Both materials are isotropic linear elastic. The parameters and the geometry of the unit-cell are shown in  Figure \ref{Fig;unit_cell}. We consider $n$ unit-cells in each direction in the computation domain $\B = [-\dfrac{L}{2},\dfrac{L}{2}] \times  [-\dfrac{L}{2},\dfrac{L}{2}]  =  [-\dfrac{n \, l}{2},\dfrac{n \, l}{2}] \times  [-\dfrac{n \, l}{2},\dfrac{n \, l}{2}]$.  

\section{Motivation and consistency check for linear elasticity} 
\label{section:linearelasticity}
In this section, we introduce the fundamentals of a homogenization procedure based on the concept of a least squares fitting of energies for identifying the unknown parameters of the effective homogenized continuum.

\subsection{Affine Dirichlet boundary condition} 
\label{section:linearelasticity:affune}
\begin{figure}[htpb]
\center
	\unitlength=1mm
	\begin{picture}(140,40)
	\put(0,0){\def\svgwidth{14 cm}{\small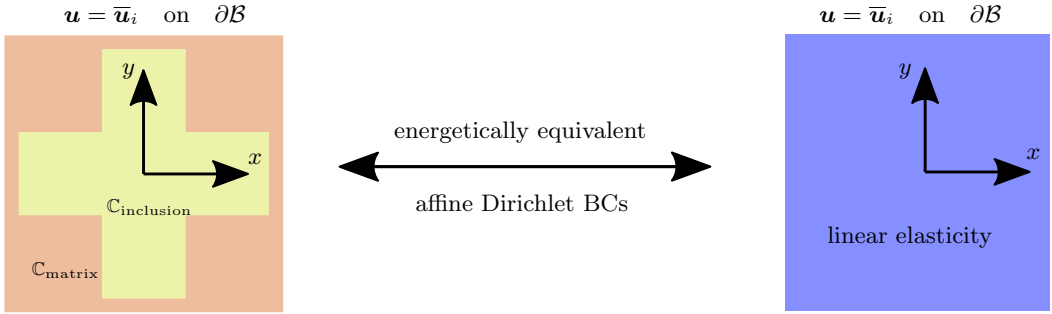}}
	\end{picture}
\caption{Illustration of the homogenization procedure to identify the parameters of an equivalent linear elastic medium under affine kinematic boundary condition.}
\label{Figure:illustration:linear_elasticity}
\end{figure}

We employ this approach to first identify the elasticity tensor of an equivalent linear elastic homogeneous medium under affine boundary conditions, as shown in Figure \ref{Figure:illustration:linear_elasticity}.  
While affine Dirichlet boundary conditions are usually not considered in standard homogenization approaches, it is nevertheless true that there exists a unique effective elasticity tensor $\mathbb{C}^\textrm{affine}$ which defines energy equivalence under all affine boundary conditions which is  typically stiffer than the effective elasticity tensor $\mathbb{C} ^\textrm{periodic}$ obtained under periodic boundary condition,  see \cite{NefEidMad:2019:ios}.  We use the generalized notation $\mathbb{C} ^\textrm{hom}$ to encompass both cases since the same algorithm is involved.

 Here, $ \Pi^\textrm{het}_i(\bu,\mathbb{C})$ is the total energy of the heterogeneous domain  under loading case $i$ while $ \Pi^\textrm{hom}_i (\bu, \mathbb{C}^\textrm{hom}) $ is the total energy of an equivalent linear elastic homogeneous continuum under loading case $i$. They read

\begin{equation}
\Pi^\textrm{het}_i (\bu, \mathbb{C}  )   =  \int_{\B}  \Bvarepsilon : \mathbb{C} (\bx) :\Bvarepsilon \, \textrm{d} V \,, \qquad   \Pi^\textrm{hom}_i (\bu, \mathbb{C}^\textrm{hom})  =  \int_{\B}   \Bvarepsilon : \mathbb{C}^\textrm{hom} :\Bvarepsilon \, \textrm{d} V \,. 
\end{equation}

The Hill-Mandel lemma postulates that the mechanical energies of the heterogeneous and equivalent homogeneous continua are equivalent.   We have the following minimization problem 

\begin{equation}
\label{Eq:minimization1}
r^2 = \min_{ \mathbb{C}^\textrm{hom}}  \sum_{i=1}^{i_\textrm{max}}  \, \lvert \lvert  \Pi^\textrm{het}_i(\bu,  \mathbb{C}  ) -\Pi^\textrm{hom}_i(\bu, \mathbb{C}^\textrm{hom}) \rvert \rvert^2  \,, 
\end{equation}

where $i$ indicates different loading cases  induced by affine Dirichlet boundary conditions enforced on the whole boundary (${\bar{\bu}_i = \bar{\Bvarepsilon}_i \cdot \bx \quad \textrm{on} \quad \partial \B } \quad$ with $ \quad \bar{\Bvarepsilon}_i \in \Sym(3) \, $).
In our approach the stiffness's anisotropy class needs to be defined in advance. This may require additional considerations and careful selection of the appropriate stiffness's anisotropy class to ensure the accuracy and reliability of the results. We consider an equivalent continuum with a stiffness tensor which exhibits cubic anisotropy. This is inline with our preference for a cubic unit-cell in Figure  \ref{Fig;unit_cell}.  The stiffness matrix $\mathbb{C}^\textrm{hom}$ in Voigt notation reads

\begin{equation}
\widetilde{\mathbb{C}}^\textrm{hom}= \left[ \begin{array}{c c c}
2 \mu^\textrm{hom} + \lambda^\textrm{hom} & \lambda^\textrm{hom} & 0 \\ 
\lambda^\textrm{hom} & 2 \mu^\textrm{hom} + \lambda^\textrm{hom} & 0 \\
0 & 0 & {\mu^*}^\textrm{hom}
\end{array}
\right]  \,  
\end{equation}

 and is characterized by three parameters. The total energy can be rewritten as 

\begin{equation} 
\begin{aligned}
\Pi^\textrm{hom}_i (\bu,\mathbb{C}^\textrm{hom})  &=   \int_\B  \Bvarepsilon:\mathbb{C}^\textrm{hom}:\Bvarepsilon \, \textrm{d} V =  \int_\B  \mu^\textrm{hom}  \, (\varepsilon_{11}^2+\varepsilon_{22}^2) \, + {\mu^*}^\textrm{hom} ( 2 \varepsilon_{12}^2) \,  + \frac{ \lambda^\textrm{hom}}{2} (\varepsilon_{11}+\varepsilon_{22})^2 \, \textrm{d} V \\  
&=    \mu^\textrm{hom} \underbrace{\left( \int_\B  (\varepsilon_{11}^2+\varepsilon_{22}^2) \,  \textrm{d} V \right)}_{\Pi^\textrm{hom}_{i,\mu}}        +    {\mu^*}^\textrm{hom}  \underbrace{\left( \int_\B   2 \varepsilon_{12}^2 \,    \textrm{d} V \right)}_{\Pi^\textrm{hom}_{i,\mu^*}}   + \lambda^\textrm{hom}  \underbrace{\left( \int_\B \frac{1}{2} (\varepsilon_{11}+\varepsilon_{22})^2  \textrm{d} V \right)}_{\Pi^\textrm{hom}_{i,\lambda}}      \\ 
&=   \mu^\textrm{hom} \, {\Pi^\textrm{hom}_{i,\mu}}  + {\mu^*}^\textrm{hom} \, {\Pi^\textrm{hom}_{i,\mu^*}} + \lambda^\textrm{hom} \, {\Pi^\textrm{hom}_{i,\lambda}}  \,.    
\end{aligned}
\end{equation}

Our objective is to identify the parameters $\lambda^\textrm{hom}$, $\mu^\textrm{hom}$, and ${\mu^*}^\textrm{hom}$ by solving the minimization problem in Equation \ref{Eq:minimization1}. However, initial values must be set in order to solve the $n$ boundary value problems. Consequently, an iterative approach is necessary after establishing initial values for the unknown parameters. We then seek the increments $\Delta \lambda^\textrm{hom}$, $\Delta \mu^\textrm{hom}$, and $\Delta {\mu^*}^\textrm{hom}$ such that: 

\begin{equation} 
\begin{aligned}
\Pi^\textrm{hom}_i (\bu, \mathbb{C}^\textrm{hom}+ \Delta \mathbb{C}^\textrm{hom}) =& \, (\mu^\textrm{hom}+\Delta \mu^\textrm{hom})  {\Pi^\textrm{hom}_{i,\mu}} + ({\mu^*}^\textrm{hom} + \Delta {\mu^*}^\textrm{hom}) {\Pi^\textrm{hom}_{i,\mu^*}} \\ & +  (\lambda^\textrm{hom}+\Delta \lambda^\textrm{hom})  {\Pi^\textrm{hom}_{i,\lambda}}   \\ 
=& \, \Pi^\textrm{hom}_i (\bu,\mathbb{C}^\textrm{hom}) + \Delta \mu^\textrm{hom} \, {\Pi^\textrm{hom}_{i,\mu}}  + \Delta {\mu^*}^\textrm{hom} \, {\Pi^\textrm{hom}_{i,\mu^*}}   + \Delta \lambda^\textrm{hom} \, {\Pi^\textrm{hom}_{i,\lambda}}    \,.     
\end{aligned}
\end{equation}

The total energy of the homogeneous domain can be written for $i_\textrm{max}$ loading cases as 

\begin{equation}
\left[\begin{array}{c}
\Pi^\textrm{hom}_1 (\bu,\mathbb{C}^\textrm{hom} + \Delta \mathbb{C}^\textrm{hom})  \\ 
\Pi^\textrm{hom}_2 (\bu,\mathbb{C}^\textrm{hom}+ \Delta \mathbb{C}^\textrm{hom})  \\ 
. \\
. \\
\Pi^\textrm{hom}_{i_\textrm{max}} (\bu,\mathbb{C}^\textrm{hom}+ \Delta \mathbb{C}^\textrm{hom})  \\ 
\end{array}
\right]
= 
\underbrace{\left[
\begin{array}{c c c}
 \Pi^\textrm{homo}_{1,\mu} &  \Pi^\textrm{hom}_{1,\mu^*} & \Pi^\textrm{hom}_{1,\lambda} \\[1em]  
 \Pi^\textrm{hom}_{2,\mu} &  \Pi^\textrm{hom}_{2,\mu^*} & \Pi^\textrm{hom}_{2,\lambda} \\ 
. &. &. \\ 
. &. &. \\ 
\Pi^\textrm{hom}_{i_\textrm{max},\mu} &  \Pi^\textrm{hom}_{i_\textrm{max},\mu^*} & \Pi^\textrm{hom}_{i_\textrm{max},\lambda}  \\ 
\end{array}
\right]}_{\bD}
\underbrace{\left[
\begin{array}{c}
\Delta \mu^\textrm{hom} \\
\Delta {\mu^*}^\textrm{hom} \\ 
\Delta \lambda^\textrm{hom} \\
\end{array}
\right]}_{\Delta \mathbb{C}^\textrm{hom}} + 
\underbrace{\left[\begin{array}{c}
\Pi^\textrm{hom}_1 (\bu,\mathbb{C}^\textrm{hom})  \\ 
\Pi^\textrm{hom}_2 (\bu,\mathbb{C}^\textrm{hom})  \\ 
. \\
. \\
\Pi^\textrm{hom}_{i_\textrm{max}} (\bu,\mathbb{C}^\textrm{hom})  \\ 
\end{array}
\right]}_{\bb} \,. 
\end{equation}

The minimization problem becomes 

\begin{equation}
r^2=\min_{\Delta \mathbb{C}^\textrm{hom}}  \sum_{i=1}^{i_\textrm{max}}   \,  \lvert \lvert  \Pi^\textrm{het}_i (\bu,\mathbb{C}) - \left(\Pi^\textrm{hom}_i (\bu,\mathbb{C}^\textrm{hom}) + \Delta \mu^\textrm{hom}  {\Pi^\textrm{hom}_{i,\mu}} +  \Delta {\mu^*}^\textrm{hom} {\Pi^\textrm{hom}_{i,\mu^*}}  +  \Delta \lambda^\textrm{hom}   {\Pi^\textrm{hom}_{i,\lambda}}  \right) \rvert \rvert^2 \,,
\end{equation}

and the solution of the least square problem reads 

\begin{equation}
\Delta \mathbb{C}^\textrm{hom} = ( \bD^T \cdot \bD )^{-1} \cdot \bD^T  \cdot (\ba - \bb)      \qquad \textrm{with} \quad {a_i = \Pi_i^\textrm{het}(\bu, \mathbb{C} ) } \,,
\end{equation}

where the parameters  $\mu^\textrm{hom},{\mu^*}^\textrm{hom}$ and $\lambda^\textrm{hom}$ are updated in an iterative procedure until the error $r^2$ converges to a constant value, optimally zero. The algorithm is explained in Algorithm  \ref{alg:flow}. For this algorithm a standard T2 finite element with quadratic shape functions is implemented.  

\begin{algorithm}[htpb] 
  \caption{Algorithm for the minimization problem for equivalent homogeneous linear elastic continuum. }    \label{alg:flow}
  \Begin{
 \Block{ {\bf production of the reference data} (heterogeneous material)}{
   
  - {\bf inputs}: unit-cell geometry with the material parameters of the matrix and inclusion \\ 
  
  - define ${i_\textrm{max}}$ affine deformation modes $\bar{\Bvarepsilon}_i$ for $i=1,....,{i_\textrm{max}}$  \\ 
  
  - \lIf {$\textrm{affine BCs}$} {${\bu = \bar{\bu} = \bar{\Bvarepsilon}_i \cdot \bx \quad \textrm{on} \quad \partial \B }$}
  -  \lIf {periodic BCs} {${\bar{\bu} = \bar{\Bvarepsilon}_i \cdot \bx \quad \textrm{on} \quad  \B }$} 
  
  - solve $i_\textrm{max}$ boundary value problems of the heterogeneous material \\ 
  
  - calculate the vector  $\ba^T = [\Pi_1^\textrm{het}, \Pi_2^\textrm{het}, .. , \Pi^\textrm{het}_{i_\textrm{max}}] $  \\ 
	}   

 \vspace{1em}

    \Block{ {\bf defining the unknown } (homogeneous material)}{  
   
  - set initial values for the parameters $\lambda^\textrm{hom},\mu^\textrm{hom},{\mu^*}^\textrm{hom}$ \\
  
\Repeat{ $r^2 < \textrm{tol} $}{
  
  - apply the deformation modes ($\bar{\Bvarepsilon}_i$ for $i=1,....,{i_\textrm{max}}$)  \\ 
  
  - \lIf {$\textrm{affine BCs}$} {${\bu = \bar{\bu} = \bar{\Bvarepsilon}_i \cdot \bx \quad \textrm{on} \quad \partial \B }$}
  -  \lIf {periodic BCs} {${\bar{\bu} = \bar{\Bvarepsilon}_i \cdot \bx \quad \textrm{on} \quad  \B }$} 
  
 - solve $i_\textrm{max}$ boundary value problems of the equivalent homogeneous medium \\

 -  calculate the vector  $\bb^T = [\Pi_1^\textrm{hom}, \Pi_2^\textrm{hom}, .. , \Pi^\textrm{hom}_{i_\textrm{max}}] $ 
  \\ 
  
 -  calculate the derivative  matrix $\bD$; $\bD_i^T = [ \Pi^\textrm{homo}_{i,\mu} ,  \Pi^\textrm{hom}_{i,\mu^*}, \Pi^\textrm{hom}_{i,\lambda} ]$ for $i=1,..,i_\textrm{max}$ \\ 
  
- solve $ [
\Delta \mu^\textrm{hom} ,
\Delta {\mu^*}^\textrm{hom}, 
\Delta \lambda^\textrm{hom}  ]^T  = ( \bD^T \cdot \bD )^{-1} \cdot \bD^T  \cdot (\ba - \bb)$    \\  

- $ \mu^\textrm{hom} \leftarrow \mu^\textrm{hom} + \Delta \mu^\textrm{hom} , \quad {\mu^*}^\textrm{hom} \leftarrow {\mu^*}^\textrm{hom} + \Delta {\mu^*}^\textrm{hom}, \quad \lambda^\textrm{hom} \leftarrow \lambda^\textrm{hom} + \Delta \lambda^\textrm{hom}$  \\ 

- calculate the current error $r^2$  

} 
 }
 
  \vspace{1em}
 
  - The parameters $\mu^\textrm{hom},{\mu^*}^\textrm{hom},\lambda^\textrm{hom}$  are known 
  
  }
\end{algorithm}

 The implemented algorithm leads for any ${i_\textrm{max}} \ge 3$ and one single unit-cell ($n=1$) to the solution $\mu^\textrm{hom} = 6.251$ \, kN/mm$^2$, ${\mu^*}^\textrm{hom} =  \, 8.337$ kN/mm$^2$ and  $\lambda^\textrm{hom} = 4.379$ \, kN/mm$^2$ which meet the one in \cite{NefEidMad:2019:ios} within one iteration. This is expected because the system is linear. The error vanishes which means the fitting delivers actually the unique solution.  Note that it is enough to set ${i_\textrm{max}}=3$.  However, stretching in $x$ and $y$ directions leads to the same equation (energetically equivalent) for our geometry which means that the three classical affine deformation modes from the classical homogenization theory (one shear mode and two stretching modes) are not a valid choice and a third arbitrary one has to be considered.

 \begin{figure}[htpb] 
  \centering
  \fbox{
 \parbox{\textwidth}{ 
We consider four random modes $i_\textrm{max} = 4$ applied on a single unit-cell $n=1$  \\

$\bar \varepsilon_1 = \left( \begin{array}{c c } 
-0.02 & 0.03 \\
0.03 & 0.01 
\end{array}  \right) \, ,  \bar \varepsilon_2 = \left( \begin{array}{c c } 
0.03 & -0.01 \\
-0.01 & 0.05 
\end{array}  \right)  \, ,  \bar \varepsilon_3 = \left( \begin{array}{c c } 
0.01 & 0.01 \\
0.01 & -0.01 
\end{array}  \right)  \, ,  \bar \varepsilon_4 = \left( \begin{array}{c c } 
0.01 & 0 \\
0 & 0.02 
\end{array}  \right)  $ \\ \\

We assume initial values:  $ \mu^\textrm{hom} = 26.32 $ kN/mm$^2$ , ${\mu^*}^\textrm{hom} =  \, 26.32$ kN/mm$^2$ and $\lambda^\textrm{hom} = 51.08 $ kN/mm$^2$. \\

The algorithm delivers:  \\

\begin{tabular}{ccccc}
            \toprule
            iteration & $\mu^\textrm{hom}$ [kN/mm$^2$]  & ${\mu^*}^\textrm{hom}$ [kN/mm$^2$] & $\lambda^\textrm{hom}$ [kN/mm$^2$] & $r^2$ [(kN$\cdot$mm)$^2$] \\
            \midrule
               0   & $26.32 $ & $26.32 $ & $ 51.08 $ & $ 0.05198 $   \\ 
               \midrule
               1  & $6.251 $ & $8.337 $   & $ 4.379 $ & $ 2.374 \, 10^{-26} $\\     
                \midrule
               2  & $6.251 $ & $8.337 $   & $ 4.379 $ & $ 1.105 \, 10^{-28} $\\   
            \bottomrule
\end{tabular} \\ \\

 final result: $\mu^\textrm{hom} = 6.251$ \, kN/mm$^2$ , ${\mu^*}^\textrm{hom} =  \, 8.337$ kN/mm$^2$  and  $ \lambda^\textrm{hom} = 4.379$ \, kN/mm$^2$.

}}
  \caption{Results for the parameter identification algorithm for linear elasticity under affine boundary conditions.   }
    \label{fig:results}
\end{figure}

\subsection{Periodic boundary condition} 
\label{sec:PBCs}
 
In classical homogenization theory, periodicity is the choice to define the effective properties because of the significant difference in length scales between the macroscopic and microscopic problems. In essence, the micro-problem is significantly smaller than the macro-problem, thereby maintaining a clear scale separation.  The microscopic strain is decomposed into a constant macroscopic part $\bar\Bvarepsilon$ and a fluctuation  part  $\widehat\Bvarepsilon$

\begin{equation}
\Bvarepsilon = \bar\Bvarepsilon + \widehat\Bvarepsilon \qquad \textrm{and} \qquad \bu =  \bar\bu + \widehat\bu = \bar\Bvarepsilon \cdot \bx + \widehat\bu \,. 
\end{equation}

The partial derivatives of the energy with respect to the unknown parameters are computed, taking into consideration that the integral of fluctuation part of strain over the domain equals to zero, i.e. $\int_B \widehat\Bvarepsilon \, \textrm{d} V = \bzero$, 

\begin{equation}
\begin{aligned}
\Pi^\textrm{hom}_{i,\mu} &= \int_B (\varepsilon_{11}^{\,2} +  \varepsilon_{22}^{\,2}) \,\textrm{d} V  =  V_\B (\bar\varepsilon_{11}^2 + \bar\varepsilon_{22}^2) + \underbrace{\int_B 2 (\bar\varepsilon_{11}\widehat\varepsilon_{11} + \bar\varepsilon_{22}\widehat\varepsilon_{22}) \,\textrm{d} V}_{=\,0}  + \int_B  (\widehat\varepsilon_{11}^{\,2} + \widehat\varepsilon_{22}^{\,2}) \,\textrm{d} V    \,,  \\
\Pi^\textrm{hom}_{i,\mu^*} &= \int_\B   2 \varepsilon_{12}^2 \,    \textrm{d} V = 2 V_\B \bar\varepsilon_{12}     +  \underbrace{\int_B 4 \bar\varepsilon_{12} \widehat\varepsilon_{12}  \,\textrm{d} V}_{=\,0}    +  \int_B 2 \widehat\varepsilon_{12}  \,\textrm{d} V \,. \\
\Pi^\textrm{hom}_{i,\lambda} &= \int_B \frac{1}{2} (\varepsilon_{11} + \varepsilon_{22})^2 \,\textrm{d} V   \\
&= \frac{V_\B}{2} (\bar\varepsilon_{11} + \bar\varepsilon_{22})^2  + \underbrace{\int_B (\bar\varepsilon_{11} + \bar\varepsilon_{22}) (\widehat\varepsilon_{11} + \widehat\varepsilon_{22}) \, \textrm{d} V}_{=\,0}  + \int_B \frac{1}{2} (\widehat\varepsilon_{11} + \widehat\varepsilon_{22})^2 \,\textrm{d} V    \,,  \\
\end{aligned}
\end{equation}

\begin{figure}[htpb]
\center
	\unitlength=1mm
	\begin{picture}(160,45)
	\put(0,0){\def\svgwidth{16 cm}{\small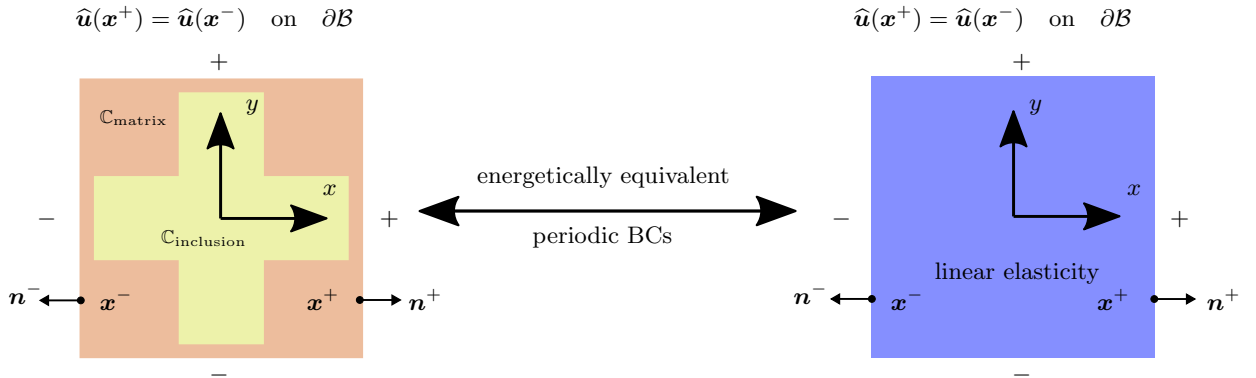}}
	\end{picture}
\caption{Illustration of periodicity condition of fluctuation on heterogeneous medium and elastic equivalent medium.}
\label{Figure:illustration:linear_elasticity_p}
\end{figure}

 The boundary $\partial \B$ is divided into two associated parts ("+","-") which satisfies $\partial \B = \partial \B^+ \cup \partial \B^- $ with outward unit normals $\bn^-$ and $\bn^+$ satisfying $\bn^- = - \bn^+$ and the periodicity is postulated as

\begin{equation}
\widehat\bu (\bx^+)   = \widehat\bu (\bx^-)     \, ,     
\end{equation} 

which is illustrated in Figure \ref{Figure:illustration:linear_elasticity_p}. The algorithm is similar to the one with affine Dirichlet boundary  conditions in Algorithm  \ref{alg:flow}, however, the $i_\textrm{max}$ deformation modes are enforced on the body, i.e. ${\bar\bu = \bar\Bvarepsilon \cdot \bx \quad \textrm{on} \quad \B}$. For this, we utilize a standard T2 finite element which discretizes the fluctuation field.

\begin{figure}[htpb]
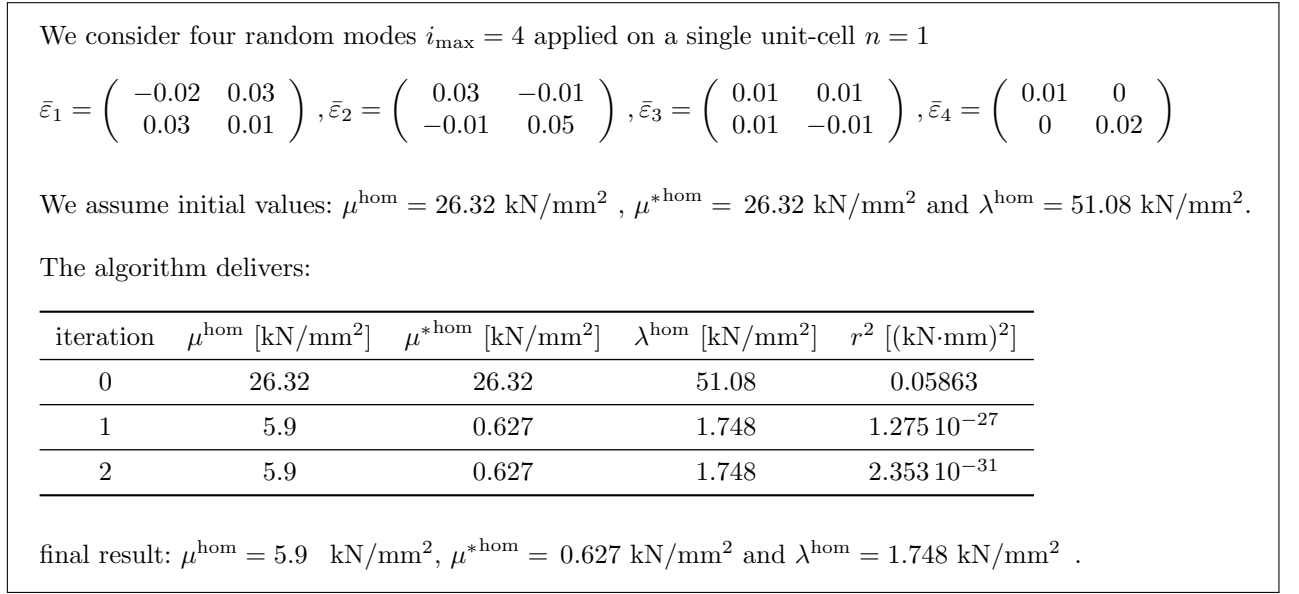
 
  \centering
  \fbox{
 \parbox{\textwidth}{ 
We consider four random modes $i_\textrm{max} = 4$ applied on a single unit-cell $n=1$  \\

$\bar \varepsilon_1 = \left( \begin{array}{c c } 
-0.02 & 0.03 \\
0.03 & 0.01 
\end{array}  \right) \, ,  \bar \varepsilon_2 = \left( \begin{array}{c c } 
0.03 & -0.01 \\
-0.01 & 0.05 
\end{array}  \right)  \, ,  \bar \varepsilon_3 = \left( \begin{array}{c c } 
0.01 & 0.01 \\
0.01 & -0.01 
\end{array}  \right)  \, ,  \bar \varepsilon_4 = \left( \begin{array}{c c } 
0.01 & 0 \\
0 & 0.02 
\end{array}  \right)  $ \\ \\

We assume initial values:  $ \mu^\textrm{hom} = 26.32 $ kN/mm$^2$ , ${\mu^*}^\textrm{hom} =  \, 26.32$ kN/mm$^2$ and $\lambda^\textrm{hom} = 51.08 $ kN/mm$^2$. \\ 

The algorithm delivers:  \\

\begin{tabular}{ccccc}
            \toprule
            iteration &  $\mu^\textrm{hom}$ [kN/mm$^2$]  & ${\mu^*}^\textrm{hom}$ [kN/mm$^2$] & $\lambda^\textrm{hom}$ [kN/mm$^2$] & $r^2$ [(kN$\cdot$mm)$^2$] \\
            \midrule
               0  & $26.32 $ & $26.32 $  & $ 51.08 $ & $ 0.05863 $   \\ 
               \midrule
               1  & $5.9 $ & $0.627 $  & $ 1.748 $ & $ 1.275 \, 10^{-27} $\\     
                \midrule
               2  & $5.9 $ & $0.627 $   & $ 1.748 $ & $ 2.353 \, 10^{-31} $\\    
            \bottomrule
\end{tabular} \\ \\

 final result:   $\mu^\textrm{hom} = 5.9$ \, kN/mm$^2$, ${\mu^*}^\textrm{hom} =  \, 0.627$ kN/mm$^2$  and $\lambda^\textrm{hom} = 1.748$ kN/mm$^2$ \,.

}}
  \caption{Results for the parameter identification algorithm for linear elasticity under periodic boundary conditions. It is observed, as expected, that $\mathbb{C}^\textrm{affine}$ is stiffer than $\mathbb{C}^\textrm{periodic}$.    }
    \label{fig:results_p}
\end{figure}

The results of the implemented algorithm are shown in  Figure \ref{fig:results_p}. For any selection of three deformation modes or more $i_\textrm{max} \ge 3$, we obtain the exact solution in one iteration which reads $\mu^\textrm{hom} = 5.9$ \, kN/mm$^2$, ${\mu^*}^\textrm{hom} =  \, 0.627$ kN/mm$^2$ and  $\lambda^\textrm{hom} = 1.748$ \, kN/mm$^2$ and meets the solution in \cite{SarSchSchNef:2023:oti,NefEidMad:2019:ios}. Similar to the case of affine boundary conditions, the classical three deformation modes (stretching in $x$, stretching in $y$ and shearing) are not a valid option because stretching in $x$ or $y$ leads to the same equation (energetically equivalent).  The result of this analysis serves as the limit case for the relaxed micromorphic model for very large specimens  $n \rightarrow \infty$, i.e. linear elasticity with elasticity tensor $\Cmacro$.  Enforcing vanishing fluctuations on the boundary restores  the same results  of the affine Dirichlet boundary condition, of course.

\section{Computational approach to identify the material parameters for RMM} 
\label{section:RMM}

The macro elasticity tensor $\Cmacro$ corresponds to the case when $n \rightarrow \infty$, where the macro homogeneous response is expected. A unit-cell with periodic boundary conditions should be used to identify $\Cmacro$, as described in \cite{ZohWri:2005:ait}.  Our analysis in Section \ref{sec:PBCs}
ended  with the following macroscopic stiffness exhibiting cubic symmetry 

\begin{equation}
\label{eq:macro}
\CmacroV = \left[ \begin{array}{c c c}
2 \mu_\textrm{macro} + \lambda_\textrm{macro} & \lambda_\textrm{macro} & 0 \\ 
\lambda_\textrm{macro} & 2 \mu_\textrm{macro} + \lambda_\textrm{macro} & 0 \\
0 & 0 & \mu^*_\textrm{macro} 
\end{array} \right]    \,,   \begin{array}{c c c c}
\lambda_\textrm{macro} & = &  \, 1.748 \, & \textrm{ kN/mm$^2$}\\ 
\mu_\textrm{macro} & = & 5.9 \, & \textrm{ kN/mm$^2$}\\ 
\mu^*_\textrm{macro} & =  & \, 0.627 \, & \textrm{ kN/mm$^2$}\\
\end{array}\,,
\end{equation}

where the cubic stiffness tensor is defined by three independent parameters.  Based on the extended Neumann's principle in \cite{NefEidMad:2019:ios}, the stiffness tensors $\Ce$ and $\Cmicro$ must contain the maximal invariance
group of the periodic microstructure ( i.e. $\Cmacro$). Thus, we assume 

\begin{equation}
\label{Eq:tensors:voigtnotation}
 \CmicroV= \left[ \begin{array}{c c c}
2 \mu_\textrm{micro} + \lambda_\textrm{micro} & \lambda_\textrm{micro} & 0 \\ 
\lambda_\textrm{micro} & 2 \mu_\textrm{micro} + \lambda_\textrm{micro} & 0 \\
0 & 0 & \mu^*_\textrm{micro}
\end{array}
\right]  , \, \quad
 \CeV= \left[ \begin{array}{c c c}
2 \mu_\textrm{e} + \lambda_\textrm{e} & \lambda_\textrm{e} & 0 \\ 
\lambda_\textrm{e} & 2 \mu_\textrm{e} + \lambda_\textrm{e} & 0 \\
0 & 0 & \mu^*_\textrm{e}
\end{array}
\right]  \, ,
\end{equation}

and the elastic energy density after considering $\Cc = \bzero$  (a choice) and $\IL = \II$ (plane strain) becomes

\begin{equation}
\begin{aligned}
W\left(\nabla \bu,\Bdis,\Curl \Bdis \right) =& \,  \mu_\textrm{e} \, \left((u_{1,1} - \dis_{11})^2  +  (u_{2,2} - \dis_{22})^2 \right)   
+ \frac{\mu^*_\textrm{e}}{2} (u_{1,2} + u_{2,1} - \dis_{12} - \dis_{21})^2   \\ 
&+  \frac{ \lambda_\textrm{e}}{2} \, (u_{1,1} + u_{2,2} - \dis_{11} - \dis_{22})^2 
 + \mu_\textrm{micro} ( \dis_{11}^2  + \dis_{22}^2 ) \\
&+ \frac{\mu^*_\textrm{micro}  }{2} (\dis_{12}+\dis_{21})^2 
 +  \frac{\lambda_\textrm{micro} }{2} (\dis_{11} + \dis_{22})^2 \\ 
&+  \frac{\mu {\Lc}^2 }{2 n^2} \left( (\Curl \Bdis)_{13}^2 + (\Curl \Bdis)_{23}^2\right) \,. 
\end{aligned}
\end{equation}

Note that assuming $\Cc = \bzero$ is a valid option, when the consistent boundary condition is applied on the  boundary, \cite{SarSchSchNef:2023:seo}. Moreover, the curvature for the 2D case is isotropic because $\Curl \Bdis$ is reduced to a vector out of the plane \cite{DagBarGhiEidNefMad:2020:edo}, see Equation \ref{Eq:Curl_2D}. Therefore, the curvature will be controlled by only one parameter with assuming that $\IL = \II$ is the fourth order identity tensor. The Reuss-like homogenization relation in Equation \ref{eq:reuss_like} taking into consideration the  relations in Equation \ref{Eq:tensors:voigtnotation} leads to  

\begin{equation}
\mu_\textrm{e} = \frac{\mu_\textrm{micro} \, \mu_\textrm{macro}}{\mu_\textrm{micro} - \mu_\textrm{macro}} , \quad \mu_\textrm{e}^* = \frac{\mu_\textrm{micro}^* \, \mu_\textrm{macro}^*}{\mu_\textrm{micro}^* - \mu_\textrm{macro}^*} , \quad \lambda_\textrm{e} + \mu_\textrm{e} = \frac{ (\lambda_\textrm{micro} + \mu_\textrm{micro}) \,  (\lambda_\textrm{macro} + \mu_\textrm{macro}) }{(\lambda_\textrm{micro} + \mu_\textrm{micro}) -  (\lambda_\textrm{macro} + \mu_\textrm{macro}) } \,. 
\end{equation}

\begin{figure}[htpb]
\center
	\unitlength=1mm
	\begin{picture}(160,35)
	\put(0,0){\def\svgwidth{16 cm}{\small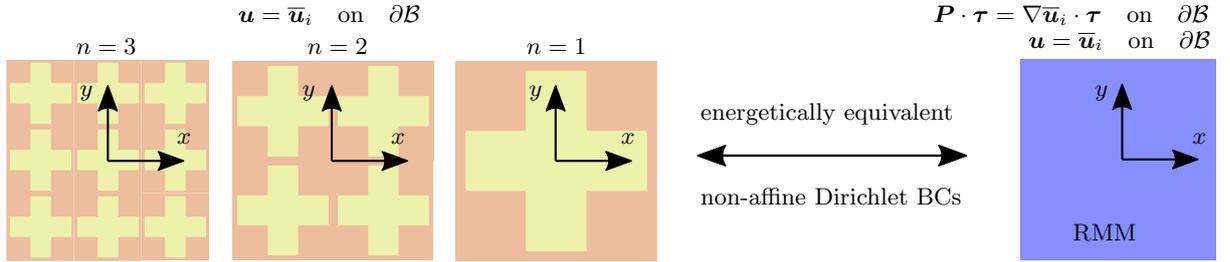}}
	\end{picture}
\caption{Illustration of the minimization problem to obtain the material parameters in the RMM. The consistent coupling condition is enforced on the whole boundary. }
\label{Figure:illustration:RMM}
\end{figure}

The presented minimization problem, illustrated in Figure \ref{Figure:illustration:RMM}, can be described as 

\begin{equation}
\label{Eq:minimization2} 
r^2 = \min_{\mu_\textrm{micro}, \, \mu^*_\textrm{micro}, \, \lambda_\textrm{micro}, \, \mu \Lc^2}  \sum_{n=1}^{n_\textrm{max}}  \sum_{i=1}^{i_\textrm{max}} \lvert \lvert  \Pi^\textrm{het}_{i \times n} (\bu) -\Pi_{i \times n} (\bu,\Bdis)  \rvert \rvert^2 \,,     
\end{equation}

where $i=1,....,i_\textrm{max}$ represents loading cases on $n \times n$ unit-cells for $n=1,...,n_\textrm{max}$. Consequently, $j_\textrm{max} = i_\textrm{max}  \, n_\textrm{max}$ reference data points need to be obtained for the heterogeneous material.

 In order to solve the minimization problem in Equation \ref{Eq:minimization2}, we define the following derivatives 

\begin{equation}
\begin{aligned}
\Pi_{{i \times n},\mu_\textrm{micro}}   &=  \frac{\partial \Pi_{i \times n}}{\partial \mu_\textrm{micro}} = \int_\B \frac{\partial W}{\partial \mu_\textrm{micro}} \, \textrm{d} V \,, \\
\Pi_{{i \times n},\mu_\textrm{micro}^*}   &=  \frac{\partial \Pi_{i \times n}}{\partial \mu_\textrm{micro}^*} = \int_\B \frac{\partial W}{\partial \mu_\textrm{micro}^*} \, \textrm{d} V \,, \\
\Pi_{{i \times n},\lambda_\textrm{micro}}   &=  \frac{\partial \Pi_{i \times n}}{\partial \lambda_\textrm{micro}} = \int_\B \frac{\partial W}{\partial \lambda_\textrm{micro}} \, \textrm{d} V \,, \\
 \Pi_{{i \times n},\mu \Lc^2}   &=  \frac{\partial \Pi_{i \times n}}{\partial \mu \Lc^2} = \int_\B \frac{\partial W}{\partial \mu \Lc^2}  \, \textrm{d} V \,, 
\end{aligned}
\end{equation}

which can not be evaluated analytically. Therefore, they will be defined numerically as

\begin{equation}
\frac{\partial \Pi}{\partial \bullet} = \frac{\Pi(\bullet + \epsilon) - \Pi(\bullet)}{ \epsilon} 
\end{equation}

where  the scalar $\epsilon$ has to be small. We reformulate the minimization problem in Equation \ref{Eq:minimization2} in terms of the increments of the unknowns instead, leading to

\begin{equation}
\label{Eq:error}
\begin{aligned}
r^2=& \min_{
\Delta \mu_\textrm{micro},
\Delta \mu^*_\textrm{micro},
\Delta \lambda_\textrm{micro},
\Delta   \mu \Lc^2  }    \sum_{n=1}^{n_\textrm{max}}  \sum_{i=1}^{i_\textrm{max}}  \, \lvert \lvert  \Pi^\textrm{het}_{i \times n} - \\ 
& \left( \Pi_{i \times n}  +    \frac{\partial \Pi_{i \times n}}{\partial \mu_\textrm{micro}}   \Delta \mu_\textrm{micro}+    \frac{\partial \Pi_{i \times n}}{\partial \mu_\textrm{micro}^*}   \Delta \mu_\textrm{micro}^* +  \frac{\partial \Pi_{i \times n}}{\partial \lambda_\textrm{micro}}  \Delta \lambda_\textrm{micro}    +
  \frac{\partial \Pi_{i \times n}}{\partial \mu \Lc^2}   \Delta   \mu \Lc^2  
  \right)  \rvert \rvert^2 \,. 
  \end{aligned}
\end{equation}

Hence, we obtain an optimization problem where the unknowns have to be updated in an iterative procedure.  The solution of the minimization problem at the current unknowns $(\mu_\textrm{micro},\mu^*_\textrm{micro},\lambda_\textrm{micro},\mu  \Lc^2)$  leads to the following vector  

\begin{equation}
\BLambda = ( \bD^T \cdot \bD )^{-1} \cdot \bD^T  \cdot (\ba - \bb)    \,, 
\end{equation}

with

\begin{equation}
\ba = \left[\begin{array}{c}
\Pi^\textrm{het}_{1 \times 1}  \\ 
\Pi^\textrm{het}_{1 \times 2}  \\ 
. \\
. \\
\Pi^\textrm{het}_{i_\textrm{max} \times n_\textrm{max}}  \\ 
\end{array}
\right] \,, \quad 
\bb = \left[\begin{array}{c}
\Pi_{1 \times 1}  \\ 
\Pi_{1 \times 2}  \\ 
. \\
. \\
\Pi_{i_\textrm{max} \times n_\textrm{max}}  \\ 
\end{array}
\right] \,, \quad  
\BLambda = 
\left[
\begin{array}{c}
\Delta \mu_\textrm{micro} \\
\Delta \mu^*_\textrm{micro} \\
\Delta \lambda_\textrm{micro} \\
\Delta \mu \Lc^2  
\end{array}
\right] \,, 
\end{equation}
and
\begin{equation}
\bD =  
\left[
\begin{array}{c}
\bD_{1 \times 1}^T \\[0.6 em] 
\bD_{1 \times 2}^T \\[0.6 em] 
. \\
. \\[0.6 em] 
\bD_{i_\textrm{max} \times n_\textrm{max}}^T 
\end{array}
\right] = 
\left[
\begin{array}{c c c c}
  \frac{\partial \Pi_{1 \times 1}}{\partial \mu_\textrm{micro}}  &   \frac{\partial \Pi_{1 \times 1}}{\partial \mu_\textrm{micro}^*} & \frac{\partial \Pi_{1 \times 1}}{\partial \lambda_\textrm{micro}} &  \frac{\partial \Pi_{1 \times 1}}{\partial \mu  \Lc^2}  \\[1 em] 
  \frac{\partial \Pi_{1 \times 2}}{\partial \mu_\textrm{micro}}  &   \frac{\partial \Pi_{1 \times 2}}{\partial \mu_\textrm{micro}^*} & \frac{\partial \Pi_{1 \times 2}}{\partial \lambda_\textrm{micro}} &  \frac{\partial \Pi_{1 \times 2}}{\partial \mu  \Lc^2}  \\
. & . & . & . \\
. & . & . & .  \\ 
  \frac{\partial \Pi_{i_\textrm{max} \times n_\textrm{max}}}{\partial \mu_\textrm{micro}}  &   \frac{\partial \Pi_{i_\textrm{max} \times n_\textrm{max}}}{\partial \mu_\textrm{micro}^*} & \frac{\partial \Pi_{i_\textrm{max} \times n_\textrm{max}}}{\partial \lambda_\textrm{micro}} &  \frac{\partial \Pi_{i_\textrm{max} \times n_\textrm{max}}}{\partial \mu  \Lc^2} \\
\end{array}
\right] \,, 
\end{equation}

where the vector $\BLambda$ represents a preferred direction at the current position.  The new position has to be updated as 

\begin{equation}
\left[ \begin{array}{c}
\mu_\textrm{micro} \\ 
\mu^*_\textrm{micro} \\
\lambda_\textrm{micro} \\  
\mu  \Lc^2
 \end{array} \right]_\textrm{new} =  \left[ \begin{array}{c}
\mu_\textrm{micro} \\ 
\mu^*_\textrm{micro} \\ 
\lambda_\textrm{micro} \\ 
\mu  \Lc^2
 \end{array} \right] + \beta \BLambda
\end{equation}
where $\beta$ is the distance along the direction $\BLambda$.  It can be identified by a line search procedure. An exact procedure can be implemented to optimize the choice of $\beta$ with the criterion

\begin{equation}
r^2 = \min_{\beta}  \sum_{n=1}^{n_\textrm{max}}  \sum_{i=1}^{i_\textrm{max}} \lvert \lvert  \Pi^\textrm{het}_{i \times n} -\Pi_{i \times n} (  \mu_\textrm{micro} +  \beta \Lambda_1 ,  \mu^*_\textrm{micro} +  \beta  \Lambda_2, \lambda_\textrm{micro} + \beta  \Lambda_3, \mu \Lc^2 +   \beta \Lambda_4)  \rvert \rvert^2 \,. 
\end{equation}

However, we need to keep $\Cmicro$ stiffer than $\Cmacro$, i.e. $\Cmicro - \Cmacro$ must be positive definite, see  \cite{NefEidMad:2019:ios}, and $\Lc$ must be strictly positive which yields then a maximum distance $\beta_\textrm{max}$ along the preferred  direction which satisfies 

\begin{equation}
\begin{aligned}
 \quad  \mu_\textrm{micro} +  \beta_\textrm{max} \Lambda_1 &> \mu_\textrm{macro}  \,, \\
  \quad \mu^*_\textrm{micro} +  \beta_\textrm{max} \Lambda_2 &> \mu^*_\textrm{macro}\,, \\
 \lambda_\textrm{micro} + \beta_\textrm{max} \Lambda_3  +  \mu_\textrm{micro} +  \beta_\textrm{max} \Lambda_1  &> \lambda_\textrm{macro} + \mu_\textrm{macro}  \,, \\  
   \mu \Lc^2 +   \beta_\textrm{max} \Lambda_4 &>  0 \,.
 \end{aligned}
\end{equation}

For our implementation, inexact identification of $\beta$ is implemented by evaluating the function at multiple points in the domain $ \beta \in \left[0,\min(1,\beta_\textrm{max})\right]$ along the preferred  direction $\Lambda$ and we choose the one with the least error $r^2$.

The choice of the boundary conditions plays a crucial role in the homogenization theory. 
We select Dirichlet boundary conditions that encompass both affine and non-affine parts on the entire boundary 

\begin{equation}
\label{eq:bcomponents}
{\bar{\bu}_i =   \bB_i \cdot \bx + \bC_i \cdot \bx \otimes \bx  \quad \textrm{on} \quad \partial \B }
\end{equation}

with 

\begin{equation}
\label{eq:ccomponents}
 \bB_i = \left[ \begin{array}{c c}
B_{11} &  B_{12} \\ 
B_{21} &  B_{22} 
 \end{array} \right]  \quad  \textrm{with}  \quad B_{jk} =  \,\, \textrm{random}[-0.05,0.05]
\end{equation}

and 

\begin{equation}
 \bC_i \cdot \bx \otimes \bx = \left[ \begin{array}{c c c}
C_{111} &  C_{112} &  C_{122}  \\ 
C_{211} &  C_{212} &  C_{222} 
\end{array} \right]  \cdot \left[ \begin{array}{c} 
x^2 \\ 
xy \\
y^2   
\end{array} \right]  
  \quad  \textrm{with}  \quad C_{jkl} =  \,\, \textrm{random}[-0.05,0.05] \,. 
\end{equation}

Certainly, the incorporation of non-affine boundary conditions can be expected, in analogy to the classical theories of higher-order homogenization. Nevertheless, size-effects have been demonstrated for both affine and non-affine loading. The implemented algorithm can be found in Algorithm \ref{alg:flow2} which is searching for an optimized quantity $\mu \Lc^2$. To illustrate  the order of the characteristic length $\Lc$ with respect to the size of the unit-cell, we give a definition for the shear modulus $\mu$ by using the isotropic elastic moduli closest to the cubic macroscopic moduli $\Cmacro$ for the log-Euclidean norm \cite{Nor:2006:tim}. This definition is unique and independent of whether the difference in stiffness or compliance is considered. The isotropic equivalent shear modulus reads

\begin{equation}
\label{eq:mue}
\mu = \sqrt[5]{ ( \mu_\textrm{macro})^2 (\mu^*_\textrm{macro})^3}= 1.537 \, \textrm{kN/mm}^2 \,. 
\end{equation}

\begin{algorithm}[htpb] 
  \caption{Algorithm for the optimization procedure of parameter identification of the relaxed micromorphic model. } \label{alg:flow2}
  \Begin{
 \Block{ {\bf production of the reference data} (heterogeneous material)}{
   
  - define $n_\textrm{max}$ deformation modes ($i=1,....,i_\textrm{max}$) applied on the boundary \\ 
   $${\bar{\bu}_i = \bB_i \cdot \bx + \bC_i \cdot \bx \otimes \bx   \quad \textrm{on} \quad \partial \B }$$ \\
   - define number of  cluster sizes considered in analysis with $n \times n$ unit-cells; $n=1,..,n_\textrm{max}$ \\
   
   - solve $j_\textrm{max} = i_\textrm{max}  n_\textrm{max}$  boundary value problems of the heterogeneous material \\
   
  - calculate the vector $\ba^T = [\Pi^\textrm{het}_{1\times 1}, \Pi^\textrm{het}_{1\times 2}, .... , \Pi^\textrm{het}_{{i_\textrm{max} \times n_\textrm{max}} } ]$ \\ 
}
 \vspace{1em}

    \Block{   {\bf defining the unknown } (homogeneous relaxed micromorphic continnum) }{  
   
  - give initial values for the parameters $\mu_\textrm{micro},\mu^*_\textrm{micro},\lambda_\textrm{micro}, \mu \Lc^2$ \\

\Repeat{ ${\dfrac{r^2_\textrm{old}-r^2_\textrm{new}}{r^2_\textrm{old}}} < \textrm{tol}$}{
  
  - apply the deformation modes on the boundary  (${\bar{\bu}_i =   \bB_i \cdot \bx + \bC_i \cdot \bx \otimes \bx  \quad \textrm{on} \quad \partial \B }$) \\ 
  
 - calculate the vector $\bb^T = [\Pi_{1 \times 1}, \Pi_{1\times 2}, .... , \Pi_{{i_\textrm{max} \times n_\textrm{max}}} ]$  \\
 
 - calculate the  matrix $\bD$; row vectors read $ \bD_{i \times n}^T = [   \frac{\partial \Pi_{i \times n}}{\partial \mu_\textrm{micro}}  ,   \frac{\partial \Pi_{i \times n}}{\partial \mu_\textrm{micro}^*} ,   \frac{\partial \Pi_{i \times n}}{\partial \lambda_\textrm{micro}} , \frac{\partial \Pi_{i \times n}}{\partial \mu \Lc^2} \ ]$ \\
  
- solve $ \BLambda  = ( \bD^T \cdot \bD )^{-1} \cdot \bD^T  \cdot (\ba - \bb)$    \\

- define $\beta_\textrm{max} \le 1 $ which keeps $\Cmicro$ positive definite and $\Lc$ positive \\

- try multiple values of $\beta = \{ \dfrac{1}{512}, \dfrac{1}{256},  \dfrac{1}{128},  \dfrac{1}{64},  \dfrac{1}{32},  \dfrac{1}{8},  \dfrac{1}{4},  \dfrac{1}{2}, 1 \} \beta_\textrm{max} $ \\

- choose  $\beta$ which delivers the least error $r^2$ along the direction $\BLambda$ \\ 

- $$ \mu_\textrm{micro} \leftarrow \mu_\textrm{micro} +  \beta \Lambda_1   , \qquad  \mu^*_\textrm{micro} \leftarrow \mu^*_\textrm{micro} +  \beta \Lambda_2$$

 $$\lambda_\textrm{micro} \leftarrow \lambda_\textrm{micro} + \beta \Lambda_3,    \qquad  \mu  \Lc^2 \leftarrow \mu \Lc^2 +   \beta \Lambda_4$$  \\

- calculate the current error $r_\textrm{new}^2$  to compare with the one form last iteration  $r_\textrm{old}^2$ \\
} 
 }
 
  \vspace{1em}
 
  - the parameters $\mu_\textrm{micro},\mu^*_\textrm{micro},\lambda_\textrm{micro},\mu  \Lc^2$  are known  \\
  
  }
\end{algorithm}

The results of the implemented algorithm for the parameter identification of the RMM are shown in Figure \ref{fig:results:RMM}. We utilize 40 deformation modes across three distinct sizes. Consequently, we need to solve 120 boundary value problems for the heterogeneous case before initiating the least square fitting procedure.

 \begin{figure}[htpb] 
  \centering
  \fbox{
 \parbox{\textwidth}{ 
We consider forty random modes $i_\textrm{max} = 40$  applied on three sizes $n_\textrm{max}=3$  

$$
{\bar{\bu}_i =   \bB_i \cdot \bx + \bC_i \cdot \bx \otimes \bx  \quad \textrm{on} \quad \partial \B }
$$

where components of $\bB_i$ and $\bC_i$ are randomly generated based on Equations \ref{eq:bcomponents} and \ref{eq:ccomponents} \\

known from previous analysis in Equations \ref{eq:macro} and \ref{eq:mue}: \\

$ \mu_\textrm{macro}  =  5.9 \,  \textrm{ kN/mm$^2$}  \,, \quad 
\mu^*_\textrm{macro}  =    0.627 \,  \textrm{ kN/mm$^2$} \,,  \quad 
\lambda_\textrm{macro}  =   1.748 \,  \textrm{ kN/mm$^2$}\,, \quad
\mu  =    1.537 \,  \textrm{ kN/mm$^2$}$   \\ 

set initial values:  \\

 $  \mu_\textrm{micro} = 26.32 $ kN/mm$^2, \quad  \mu_\textrm{micro}^* =  \, 26.32$ kN/mm$^2, \quad \lambda_\textrm{micro} = 51.08 $ kN/mm$^2 , \quad \Lc = 1 \, \textrm{mm}$ \\ 

The algorithm delivers:  \\

\begin{tabular}{cccccc}
            \toprule
            iteration &  $\mu_\textrm{micro}$ [kN/mm$^2$]  & $\mu^*_\textrm{micro}$ [kN/mm$^2$]  &  $\lambda_\textrm{micro}$ [kN/mm$^2$] & $\Lc$ [mm]  &  $r^2$ [(kN$\cdot$mm)$^2$] \\
            \toprule
               0  &  $26.32 $ & $26.32 $  &  $ 51.08 $ & $ 1 $ & $0.015237$  \\ 
               \midrule
               1  &  $19.66$ & $50.85$  & $ 38.46$ & $ 0.788$ & $0.00273407$\\     
                \midrule
               2  & $15.38 $ & $143.12$  &  $ 27.98 $ & $ 0.697$ & $0.0011551$\\   
               \midrule
               :  &    : & : & : & : & : \\ 
               \midrule
               6  &  $10.18 $ & $344.45$  & $ 11.29 $ & $ 0.883$ & $0.000352264$\\   
              \midrule
               7  &  $10.19 $ & $358.87$  & $ 11.3 $ & $ 0.882$ & $0.000352236$\\ 
                 \midrule  
                  :   & : & : & : & : & : \\ 
                 \midrule
               10  &  $10.19 $ & $354.77$  &  $ 11.3 $ & $ 0.882$ & $0.000352234$\\     
               \midrule
               11  &  $10.19 $ & $354.88$  & $ 11.3 $ & $ 0.882$ & $0.000352234$\\                
               \midrule
               12  &  $10.19 $ & $354.87$  &  $ 11.3 $ & $ 0.882$ & $0.000352234$\\           
            \bottomrule
\end{tabular} \\ \\

final parameters set: \\ 

$\mu_\textrm{micro} = 10.19$ kN/mm$^2$, $\mu^*_\textrm{micro} = 354.87$ kN/mm$^2$,  $\lambda_\textrm{micro} = 11.3$ kN/mm$^2$  and $\Lc  = 0.882 $ mm $= 0.882 \, L$.

}}
  \caption{Results for the parameter identification algorithm for the relaxed micromorphic model.  We refer to the obtained final parameters here as {\bf parameter set 1}. }
    \label{fig:results:RMM}
\end{figure}

The algorithm achieves the intended purpose, resulting in a final error that is much less than the error of the initial values. However, we observe that $\mu_\textrm{micro}^*$ tends towards very large values. Thus, we need to invoke the concept of the stiffest possible response, where we assume that the micro elasticity tensor cannot possibly yield a stiffer response than the homogeneous stiff matrix for any loading scenario. The concept of the stiffest possible response arises from the fact that the relaxed micromorphic model operates between two scales, ensuring the preservation of the physical interpretation of the upper bound, which cannot be logically stiffer than the stiff matrix.  This will be expressed in terms of energy norms (L\"owner order) as  

\begin{equation}
\Bvarepsilon : \Cmicro :  \Bvarepsilon  \quad  \le \quad \Bvarepsilon : \Cmatrix :  \Bvarepsilon \, \,, \quad  \forall   \Bvarepsilon \in \Sym(3) \, . 
\end{equation}

This condition can be rewritten as 

\begin{equation}
\begin{aligned}
& \left[  \begin{array}{c}
\varepsilon_{12} \\ 
\varepsilon_{22} \\
2 \varepsilon_{12}
\end{array} \right]^T  \cdot \left[ \begin{array}{c c c}
2 \mu_\textrm{micro} + \lambda_\textrm{micro} & \lambda_\textrm{micro} & 0 \\ 
\lambda_\textrm{micro} & 2 \mu_\textrm{micro} + \lambda_\textrm{micro} & 0 \\
0 & 0 & \mu^*_\textrm{micro}
\end{array}
\right]  \cdot  \left[  \begin{array}{c}
\varepsilon_{12} \\ 
\varepsilon_{22} \\
2 \varepsilon_{12}
 \end{array} \right]   \le     \\ 
& \left[  \begin{array}{c}
\varepsilon_{12} \\ 
\varepsilon_{22} \\
2 \varepsilon_{12}
\end{array} \right]^T  \cdot \left[ \begin{array}{c c c}
2 \mu_\textrm{matrix} + \lambda_\textrm{matrix} & \lambda_\textrm{matrix} & 0 \\ 
\lambda_\textrm{matrix} & 2 \mu_\textrm{matrix} + \lambda_\textrm{matrix} & 0 \\
0 & 0 & \mu^*_\textrm{matrix}
\end{array}
\right]  \cdot  \left[  \begin{array}{c}
\varepsilon_{12} \\ 
\varepsilon_{22} \\
2 \varepsilon_{12}
 \end{array} \right] \, \,, \, \forall   \left[  \begin{array}{c}
\varepsilon_{12} \\ 
\varepsilon_{22} \\
\varepsilon_{12}
 \end{array} \right] \in \R^3 \,, 
\end{aligned}
\end{equation}

and the solution reads, see \cite{NefEidMad:2019:ios},  

\begin{equation}
\label{eq:constrains}
\begin{aligned}
\mu^*_\textrm{micro} & \le  \, \mu_\textrm{matrix} \,, \\ 
\mu_\textrm{micro} & \le  \, \mu_\textrm{matrix} \,, \\
 \mu_\textrm{micro} + \lambda_\textrm{micro} & \le  \, \mu_\textrm{matrix} + \lambda_\textrm{matrix} \,. \\
\end{aligned}
\end{equation}

The first optimization delivers a shear modulus $\mu_\textrm{micro}^* = 354.87$  kN/mm$^2$ which does not meet the criterion for the stiffest response in Equation \ref{eq:constrains} since $\mu_\textrm{matrix}^* = 26.32$   kN/mm$^2$. Consequently, we refine the optimization algorithm. In each iteration, if the new values violate the constraints in Equation \ref{eq:constrains}, we project the parameter that breaks the constraints back into the admissible domain. We then repeat the current iteration, excluding this parameter. However, in the next iteration, we include all parameters again. In our algorithm, only the parameter $\mu_\textrm{micro}^*$ is attempting in each iteration to break the upper limit, leading to its projection back into the admissible domain in each iteration. The results of the modified algorithm are depicted in Figure \ref{fig:results:RMM:2}.  Obviously, the obtained parameters from the optimization procedure with upper constraints leads to a relatively larger error compared to when no constraints are considered. There is a concern, particularly with gradient-based algorithms, regarding whether it converges to a global minimum and the potential for a much better solution. To explore this, we computed the error for $11^4$ parameter sets, employing 10 divisions within the admissible domain for each parameter between the macro and matrix parameters. The set demonstrating the least error was chosen as the algorithm's starting point. The algorithm consistently converges to the same solution as shown in Figure \ref{fig:results:RMM:2} for any starting point (we have tried many other starting points).  The advantageous behavior of the relaxed micromorphic model is recognized as a two-scale model. It is bounded by two limits, each with distinct physical interpretations, and all the unknown parameters have well-defined ranges. Therefore, utilizing a gradient-based optimization procedure proves to be effective. Vice versa, employing a gradient-based optimization for the classical Eringen-Mindlin micromorphic theory, lacking an upper bound, does pose challenges due to the higher number of parameters and the uncertainty about the magnitude of these unknowns. We attempted to improve the fitting by introducing an additional skew-symmetric term in the energy, i.e.  $\mu_c \abs{\abs{\skew (\nabla \bu -\Bdis)}}^2$, representing the micro-rotation coupling, where $\mu_c$ is the Cosserat modulus. However, the Cosserat couple modulus trended towards negative values and needed to be projected back to zero delivering symmetric force stress $\Bsigma$ as before. This highlights the importance of using the consistent coupling boundary condition and meets the results of our principle investigations in  \cite{SarSchSchNef:2023:seo}.

 \begin{figure}[htpb] 
  \centering
  \fbox{
 \parbox{\textwidth}{ 
We consider forty random modes $i_\textrm{max} = 40$  applied on three sizes $n_\textrm{max}=3$  

$$
{\bar{\bu}_i =   \bB_i \cdot \bx + \bC_i \cdot \bx \otimes \bx  \quad \textrm{on} \quad \partial \B }
$$

where components of $\bB_i$ and $\bC_i$ are randomly generated based on Equations \ref{eq:bcomponents} and \ref{eq:ccomponents}\\

known from previous analysis in Equations \ref{eq:macro} and \ref{eq:mue}: \\

  $\mu_\textrm{macro}  =  5.9 \,  \textrm{ kN/mm$^2$}  \,, \quad 
\mu^*_\textrm{macro}  =   \, 0.627 \,  \textrm{ kN/mm$^2$} \,, \quad  \lambda_\textrm{macro}  =  \, 1.748 \,  \textrm{ kN/mm$^2$}\,, \quad 
\mu  =   \, 1.537 \,  \textrm{ kN/mm$^2$}  $   \\ 

set initial values:  \\

$\mu_\textrm{micro} = 26.32 $ kN/mm$^2 , \quad  \mu_\textrm{micro}^* =  \, 26.32$ kN/mm$^2$, $\quad \lambda_\textrm{micro} = 51.08 $ kN/mm$^2$,  $\quad \Lc = 1 \, \textrm{mm}$ \\ 

The algorithm delivers:  \\

\begin{tabular}{cccccc}
            \toprule
            iteration & $\mu_\textrm{micro}$ [kN/mm$^2$]  & $\mu_\textrm{micro}^*$ [kN/mm$^2$]   & $\lambda_\textrm{micro}$ [kN/mm$^2$] & $\Lc$ [mm]  &  $r^2$  [(kN$\cdot$mm)$^2$] \\
            \midrule
               0   & $26.32 $  &  $26.32 $ & $ 51.08 $ & $ 1 $ & $0.015237$  \\ 
               \midrule
               1  & $19.37 $  &  $26.32 $  & $ 19.33 $ & $ 1.039 $ & $0.0048919$  \\ 
               \midrule  
               2  & $13.73 $  & $26.32 $  & $ 10.8 $ &  $ 1.083$ & $0.00124667$\\    
                \midrule
               3   & $10.07 $  &  $26.32 $ & $ 7.83 $ &  $ 1.126$ & $0.000695682$\\  
               \midrule
               4  & $10.57 $  &  $26.32 $  & $ 8.19 $ &  $ 1.122$ & $0.000674557$\\  
               \midrule
               5  & $10.55 $  &  $26.32 $  & $ 8.22 $ &  $ 1.123$ & $0.000674545$\\  
                \midrule
               6  & $10.55 $  &  $26.32 $  & $ 8.22 $ &  $ 1.123$ & $0.000674545$\\ 
               \midrule
               7  & $10.55 $  &  $26.32 $  & $ 8.22 $ &  $ 1.123$ & $0.000674545$\\ 
            \bottomrule
\end{tabular} \\ \\

different initial values:  \\

$\mu_\textrm{micro} = 10.03 $ kN/mm$^2$ $, \quad  \mu_\textrm{micro}^* =  \, 26.32$ kN/mm$^2$, $\quad \lambda_\textrm{micro} = 6.7 $ kN/mm$^2$,   $\quad \Lc = 1.36 \, \textrm{mm} $ \\ 

The algorithm delivers:  \\

\begin{tabular}{cccccc}
            \toprule
            iteration  & $\mu_\textrm{micro}$ [kN/mm$^2$]  & $\mu_\textrm{micro}^*$ [kN/mm$^2$]   & $\lambda_\textrm{micro}$ [kN/mm$^2$] & $\Lc$ [mm] &  $r^2$ [(kN$\cdot$mm)$^2$] \\
            \midrule
               0  & $10.03 $  &  $26.32 $  & $ 6.7 $ & $ 1.36 $ & $0.00153283$  \\ 
               \midrule
               1  & $10.74$  & $26.32 $ &  $ 7.9$ & $ 1.057$ & $0.000785925$\\     
                \midrule
               2   & $10.42 $  & $26.32 $ & $ 8.17 $ &  $ 1.128$ & $0.000674967$\\    
                \midrule
               3   & $10.56 $  &  $26.32 $ & $ 8.22 $ &  $ 1.122$ & $0.000674548$\\  
               \midrule
               4  & $10.55 $  &  $26.32 $  & $ 8.22 $ &  $ 1.123$ & $0.000674545$\\  
               \midrule
               5   & $10.55 $  &  $26.32 $ & $ 8.22 $ &  $ 1.123$ & $0.000674545$\\  
            \bottomrule
\end{tabular} \\ \\

final parameters set: \\

 $\mu_\textrm{micro} =10.55$ \, kN/mm$^2$, $\mu_\textrm{micro}^* =26.32$ \, kN/mm$^2$,  $\lambda_\textrm{micro} = 8.22$ \, kN/mm$^2$ and $\Lc = 1.123 \textrm{mm} = 1.123 \, L$ 
}}
  \caption{Results for the parameter identification algorithm for the relaxed micromorphic model. Here, we impose a constraint on the micro elasticity tensor $\Cmicro$ to ensure it is not stiffer than the stiff matrix. We refer to the obtained  final parameters here as {\bf parameter set 2}.}
    \label{fig:results:RMM:2}
\end{figure}

\subsection{Comparison with a skew symmetric micro-distortion field - the Cosserat case}
The relaxed micromorphic model recovers the Cosserat model for the singular limit case ${\Cmicro \rightarrow \infty}$, the micro-distortion field $\Bdis$ must then be skew-symmetric, c.f.  \cite{AlaGanSadNasAkb:2022:cmo,BleNef:2022:ssi,GhiRizMadNef:2022:cme,KhaChiMadNef:2022:eau,JeoRamMueNef:2009:zfa,JeoNef:2010:eua,NefJeo:2009:anp,NefJeoFis:2010:sio}. The energy functional of the  relaxed micromorphic model is reduced then with setting $\bA := \skew \Bdis \in \mathfrak{so}(3)$ to  

 \begin{equation}
\begin{aligned}
W_\textrm{Cosserat} \left(\nabla \bu,\bA,\Curl \bA \right) = \, \dfrac{1}{2} (& \sym \nabla \bu : \Cmacro : \sym \nabla \bu  + (\skew \nabla \bu - \bA) : \Cc : ( \skew \nabla \bu - \bA) \\
&+   \mu \, \Lc^2 \, \textrm{Curl} \bA : \IL :\textrm{Curl} \bA ) \,. \\
\end{aligned}
 \end{equation}
 
 which turns out for the cubic anisotropic case into as
 
 \begin{equation}
\begin{aligned}
W_\textrm{Cosserat} \left(\nabla \bu,\bA,\Curl \bA \right) = \, & \mu_\textrm{macro}  \, (u_{1,1}^2  + u_{2,2}^2)  
 + \, \frac{\mu^*_\textrm{macro}}{2} (u_{1,2} + u_{2,1})^2   \\ 
 &+ \frac{\lambda_\textrm{macro}}{2} (u_{1,1} + u_{2,2})^2  +   \frac{\mu_c}{2} (u_{1,2}-u_{2,1}- 2 A_{12})^2 \\ 
 &+  \frac{\mu {\Lc}^2 }{2 n^2} \left( (\Curl \bA)^2_{13}  + (\Curl \bA)^2_{23} \right) \,. 
\end{aligned}
\end{equation}

Note that $\mu_c$ must be strictly positive for the Cosserat model to be operative otherwise the coupling of the fields ($ \bu, \bA$) vanishes. We applied the optimization procedure for the Cosserat model where only two unknown parameters ($\Lc$ and $\mu_c$)  need to be identified.  The results are shown in Figure  \ref{fig:results:Cosserat}. We explored two boundary condition scenarios for the micro-distortion field: consistent boundary conditions applied to the entire boundary, and free boundary conditions. The results obtained with consistent boundary conditions showed significantly better fitting. Consequently, we focus on the results obtained under consistent boundary conditions in the following analysis.

  \begin{figure}[htpb] 
  \centering
  \fbox{
 \parbox{\textwidth}{ 
We consider forty random modes $i_\textrm{max} = 40$  applied on three sizes $n_\textrm{max}=3$  

$$
{\bar{\bu}_i =   \bB_i \cdot \bx + \bC_i \cdot \bx \otimes \bx  \quad \textrm{on} \quad \partial \B }
$$

where components of $\bB_i$ and $\bC_i$ are randomly generated based on Equations \ref{eq:bcomponents} and \ref{eq:ccomponents}\\

known from previous analysis in Equations \ref{eq:macro} and \ref{eq:mue}: \\

  $\mu_\textrm{macro}  =  5.9 \,  \textrm{ kN/mm$^2$}  \,, \quad 
\mu^*_\textrm{macro}  =   \, 0.627 \,  \textrm{ kN/mm$^2$} \,, \quad  \lambda_\textrm{macro}  =  \, 1.748 \,  \textrm{ kN/mm$^2$}\,, \quad 
\mu  =   \, 1.537 \,  \textrm{ kN/mm}^2  $   \\ 

further assumptions: \\

$\mu_\textrm{micro} =  10000 \, \mu_\textrm{macro}  \,, \quad  \quad  \mu_\textrm{micro}^* =  10000  \, \mu_\textrm{macro}^* \,, \quad  \quad  \lambda_\textrm{micro} = 10000 \,  \lambda_\textrm{macro}  \qquad (\approx \Cmicro \rightarrow \infty)  $ \\

set initial values:  \\

 $\mu_c = 1 \, \textrm{ kN/mm}^2$   ,  $\quad \Lc = 1 \, \textrm{mm}$ \\ 

The algorithm delivers:  \\

\begin{tabular}{cc}
with consistent boundary condition & without consistent boundary condition \\ 
\begin{tabular}{cccc}
            \toprule
            iteration   & $\mu_c$   & $\Lc$  &  $r^2$\\
         & [kN/mm$^2$]  &  [mm]  &  [(kN$\cdot$mm)$^2$]  \\    
            \toprule
               $0$  & $1$ & $1$ & $0.01967$ \\ 
               \midrule
               $1$  & $4.04$ & $0.615$ & $0.00265$ \\ 
               \midrule     
               $2$  & $12.97$ & $0.59$ & $0.002517$ \\
               \midrule     
              :  & : & : & : \\               
               \midrule  
              $10$  & $1739.03$ & $0.614$ & $0.002128$ \\               
               \midrule  
              $11$  & $1408.36$ & $0.614$ & $0.002125$ \\               
               \midrule  
              : & : & : & : \\               
               \midrule
              $17$  & $452.64$ & $0.616$ & $0.00211923$ \\               
               \midrule  
              $18$  & $452.6$ & $0.616$ & $0.00211923$ \\             
            \bottomrule
\end{tabular}&  
\begin{tabular}{cccc}
            \toprule
            iteration   & $\mu_c$   & $\Lc$  &  $r^2$\\
         & [kN/mm$^2$]  &  [mm]  & [(kN$\cdot$mm)$^2$]  \\            
            \toprule
               $0$  & $1$ & $1$ & $0.0261827$ \\ 
               \midrule
               $1$  & $2.04$ & $0.711$ & $0.0258911$ \\ 
               \midrule     
               $2$  & $14.51$ & $0.51$ & $0.0257639$ \\
               \midrule     
               :  & : & : & : \\     
               \midrule  
              $10$  & $64.43 \cdot 10^{4}$ & $5.133$ & $0.0114968$ \\     
               \midrule  
              $11$  & $64.82 \cdot 10^{4}$ & $5.132$ & $0.0114965$ \\                 
              \midrule  
               : & : & : & : \\  
               \midrule  
              $19$  & $67.49 \cdot 10^{4}$ & $5.132$ & $0.0114945$ \\     
               \midrule  
              $20$  & $67.52 \cdot 10^{4}$ & $5.132$ & $0.0114945$ \\                             
            \bottomrule
\end{tabular} 
\end{tabular}
\\ \\

final parameters set (with consistent boundary condition): \\ 

 $\mu_c = 452.6$ kN/mm$^2$   and $\Lc  = 0.616   \textrm{mm} = 0.616 \, L$.  \\

}}
  \caption{Results for the parameter identification algorithm for the Cosserat model. }
    \label{fig:results:Cosserat}
\end{figure}

\subsection{Comparison with the full gradient}
Another interesting comparison can be made between the results obtained from the relaxed curvature ($\Curl \Bdis$) and the full curvature ($\nabla \Bdis$)  within the classical micromorphic model.  The energy functional of the Eringen-Mindlin micromorphic model reads

\begin{equation}
\begin{aligned}
W_\textrm{Eringen} \left(\nabla \bu,\Bdis,\nabla \Bdis \right) = \frac{1}{2} \Big( & \symb{ \nabla \bu - \Bdis} : \Ce : \symb{ \nabla \bu - \Bdis}  +   \sym \Bdis : \Cmicro: \sym \Bdis  \\
& + \skewb{ \nabla \bu - \Bdis} : \Cc : \skewb{ \nabla \bu - \Bdis} +  (\nabla \bu - \Bdis) : \mathbb{C}_\textrm{mixed} : \sym \Bdis
\\ &+  \frac{\mu \, \Lc^2}{n^2} \, \nabla \Bdis : \widetilde{\IL} :\nabla \Bdis  \Big)\,.
\end{aligned}
\end{equation}

Here, $\widetilde{\IL}$ is a sixth-order tensor and $\mathbb{C}_\textrm{mixed}$ is a fourth order tensor. For the planer cubic case, the tensor $\widetilde{\IL}$ associated with the curvature requires already the definition of 10 independent parameters \cite{DagMatLewBerDanNef:2024:}.  Thus, we have

{\footnotesize
\begin{equation}
\begin{aligned}
& \nabla \Bdis : \widetilde{\IL} :\nabla \Bdis  = \\ 
& \left( \begin{array}{c}
\dis_{11,1} \\
\dis_{12,2} \\
\dis_{22,1} \\
\dis_{21,2} \\
\dis_{22,2} \\
\dis_{21,1} \\
\dis_{11,2} \\
\dis_{12,1} \\
\end{array} \right)^T \cdot 
\left( \begin{array}{c c c c c c c c }
L_{111111} & L_{111122} & L_{111221} &  L_{111212} & 0 & 0 & 0 & 0 \\
L_{111122} & L_{122122} & L_{122221} &  L_{122212} & 0 & 0 & 0 & 0 \\
L_{111221} & L_{122221} & L_{221221} &  L_{221212} & 0 & 0 & 0 & 0 \\
L_{111212} & L_{122212} & L_{221212} &  L_{212212} & 0 & 0 & 0 & 0 \\
 0 & 0 & 0 & 0 & L_{111111} & L_{111122} & L_{111221} &  L_{111212} \\
 0 & 0 & 0 & 0 & L_{111122} & L_{122122} & L_{122221} &  L_{122212} \\
 0 & 0 & 0 & 0 & L_{111221} & L_{122221} & L_{221221} &  L_{221212}  \\
 0 & 0 & 0 & 0 & L_{111212} & L_{122212} & L_{221212} &  L_{212212} \\
\end{array} \right)
\cdot  
\left( \begin{array}{c}
\dis_{11,1} \\
\dis_{12,2} \\
\dis_{22,1} \\
\dis_{21,2} \\
\dis_{22,2} \\
\dis_{21,1} \\
\dis_{11,2} \\
\dis_{12,1} \\
\end{array} \right) \,, 
\end{aligned}
\end{equation}} \\ 

where $L_{111111} , L_{111122} , L_{111221},   L_{111212},  L_{122122} , L_{122221} ,  L_{122212} , L_{221221} ,  L_{221212} $ and $L_{212212}$ are to be determined. The total number of unknown parameters for the Eringen-Mindlin full micromorphic model equals to 14 for the 2D case (3 for $\Cmicro$, 1 for $\Cc$ and 10 for $\IL$) if we already exclude  the mixed term  ${(\nabla \bu - \Bdis) : \mathbb{C}_\textrm{mixed} : \sym \Bdis}$. Dealing with such large number of unknowns is not feasible for a gradient-based optimization at present.
 Therefore, for the sake of simplicity, we limit our analysis to the most simplified curvature formulation to reduce the number of unknown parameters, and thus the curvature will only be associated with the scalar $ \mu \, \Lc^2$, i.e. $ \widetilde{\IL} :\nabla \Bdis  = \nabla \Bdis$.  Yet, it is  clear that optimizing the coefficients of the tensor $ \widetilde{\IL}$ (with cubic symmetries) will deliver a better fit. The energy functional of the simplified micromorphic model for a cubic material turns into

\begin{equation}
\begin{aligned}
W_\textrm{Eringen}  \left(\nabla \bu,\Bdis,\nabla \Bdis \right) = \, & \,
  \mu_\textrm{e} \, \left((u_{1,1} - \dis_{11})^2  + (u_{2,2} - \dis_{22})^2 \right)   
 +   \frac{\mu^*_\textrm{e}}{2} (u_{1,2} + u_{2,1} - \dis_{12} - \dis_{21})^2  \\
 &+  \frac{\lambda_\textrm{e}}{2} (u_{1,1} + u_{2,2} - \dis_{11} - \dis_{22})^2  \\
 &+ \mu_\textrm{micro}  ( \dis_{11}^2  + \dis_{22}^2 ) + \frac{ \mu^*_\textrm{micro} }{2} (\dis_{12}+\dis_{21})^2 \\
 &+  \frac{\lambda_\textrm{micro}}{2} (\dis_{11} + \dis_{22})^2 
 + \frac{\mu_c}{2} (u_{1,2}-u_{2,1}-\dis_{12}+\dis_{21})  \\
 &+  \frac{\mu {\Lc}^2 }{2 n^2}  \abs{\abs{\nabla \Bdis}}^2 \,, 
\end{aligned}
\end{equation}

where the parameters $\mu_\textrm{micro},\mu^*_\textrm{micro},\lambda_\textrm{micro}, \mu_\textrm{c}$ and $\Lc$ need to be determined. The optimization follows the same procedure with enforcing the consistent boundary conditions on the whole boundary.  The results are illustrated in Figure     \ref{fig:results:MM}.
The optimization delivers   $\mu_\textrm{micro} = 5.959$ kN/mm$^2$  , $\mu^*_\textrm{micro} = 80.82$ kN/mm$^2$, $\lambda_\textrm{micro} = 12.06$ kN/mm$^2$ , $\mu_\textrm{c} = 1138.34$ kN/mm$^2$   and $\Lc  = 0.695 $ mm with an error {$r^2 = 0.000612487$ [(kN$\cdot$mm)$^2$]}. Notably, the simplified full micromorphic model favors an asymmetric force stress $\Bsigma$, differing from the relaxed micromorphic model. Moreover, $\Cmicro$ for the full micromorphic model is not associated with an upper limit stiffness property, which makes the obtained tensor $\Cmicro$ not being physically reflected and therefore incomparable with any measurable quantity. Therefore, a question arises as to whether the assumption of cubic symmetry for $\Cmicro$ is necessary for the full micromorphic model, given that it was originally assumed in the context of the relaxed micromorphic model and its two-scales realization.  However, the Reuss-like homogenization relation in Equation \ref{eq:reuss_like} is not valid when considering different anisotropic properties for $\Cmicro$ and $\Cmacro$. 

 \begin{figure}[htpb] 
  \centering
  \fbox{
 \parbox{\textwidth}{ 
We consider forty random modes $i_\textrm{max} = 40$  applied on three sizes $n_\textrm{max}=3$  

$$
{\bar{\bu}_i =   \bB_i \cdot \bx + \bC_i \cdot \bx \otimes \bx  \quad \textrm{on} \quad \partial \B }
$$

where components of $\bB_i$ and $\bC_i$ are randomly generated based on Equations \ref{eq:bcomponents} and \ref{eq:ccomponents}\\

known from previous analysis in Equations \ref{eq:macro} and \ref{eq:mue}: \\

  $\mu_\textrm{macro}  =  5.9 \,  \textrm{ kN/mm$^2$}  \,, \quad 
\mu^*_\textrm{macro}  =   \, 0.627 \,  \textrm{ kN/mm$^2$} \,, \quad \lambda_\textrm{macro}  =  \, 1.748 \,  \textrm{ kN/mm$^2$}\,,  
\quad 
\mu  =   \, 1.537 \,  \textrm{ kN/mm$^2$}  $   \\ 

further assumptions: \\

$\mathbb{C}_\textrm{mixed} = \bzero$ \\

set initial values:  \\

 $\mu_\textrm{micro} = 26.32 $ kN/mm$^2$ $ \,, \quad  \mu_\textrm{micro}^* =  \, 26.32$ kN/mm$^2$, $\quad \lambda_\textrm{micro} = 51.08 $ kN/mm$^2$,  $\quad \Lc = 1 \, \textrm{mm}$ \\ 

The algorithm delivers:  \\

\begin{tabular}{ccccccc}
            \toprule
            iteration & $\mu_\textrm{micro}$   & $\mu^*_\textrm{micro}$  &  $\lambda_\textrm{micro}$ &  $\mu_c$ & $\Lc$ &  $r^2$\\
                &  [kN/mm$^2$] & [kN/mm$^2$]  &  [kN/mm$^2$] &  [kN/mm$^2$]  & [mm]  & [(kN$\cdot$mm)$^2$]  \\         
            \toprule
               0  & $26.32 $ & $26.32 $  & $ 51.08 $ & $ 1 $& $ 1 $ & $0.029915$  \\        
            \toprule
               1  & $8.48 $ & $46.8 $  & $ 6.14 $ & $ 4.27 $& $ 0.7 $ & $0.0023058$  \\ 
            \toprule
               2  & $8.33 $ & $82.89 $ & $ 9.32 $  & $ 11.54 $& $ 0.621 $ & $0.0018288$  \\  
            \toprule
               :  & :  & : & :  & : & : & :  \\   
             \toprule    
               11   & $5.959 $ & $85.5 $ & $ 12.28 $ & $ 1057.11 $& $ 0.688 $ & $0.000612789$  \\  
            \toprule          
             12 & $5.959 $ & $82.99 $  & $ 12.16 $  & $ 1070.65 $& $ 0.692 $ & $0.000612667$  \\  
            \toprule     
             :  & :  & : & :  & : & : & :  \\   
             \toprule          
           25 & $5.959 $ & $80.83 $  & $ 12.06 $  & $ 1138.32 $& $ 0.695 $ & $0.000612487$  \\  
            \toprule        
           26 & $5.959 $ & $80.82 $  & $ 12.06 $  & $ 1138.34 $& $ 0.695 $ & $0.000612487$  \\  
            \bottomrule
\end{tabular} \\ \\

final parameters set: \\ 

 $\mu_\textrm{micro} = 5.959$ kN/mm$^2$  , $\quad \mu^*_\textrm{micro} = 80.82$ kN/mm$^2$,  $\quad \lambda_\textrm{micro} = 12.06$ kN/mm$^2$,   \\
 
  $\mu_c = 1138.34$ kN/mm$^2$,  $\quad \Lc  = 0.695 \, \textrm{mm} = 0.695 \, L$ .

}}
  \caption{Results for the parameter identification algorithm for the simplified full micromorphic model. }
    \label{fig:results:MM}
\end{figure} 

We can improve the fitting by introducing a simplified isotropic curvature, which is characterized by three independent parameters ($\alpha_1, \alpha_2, \alpha_3$). This isotropic curvature in the full micromorphic model has the form \cite{RizHueMadNef:2021:aso4}

\begin{equation}
\begin{aligned}
\nabla \Bdis : \widetilde{\IL} :\nabla \Bdis  =  \sum_{i=1}^2  \left( \alpha_1  \abs{\abs{\dev \sym \Bdis_{,i}}}^2 +  \alpha_2 \abs{\abs{\skew \Bdis_{,i}}}^2 + \frac{2}{9}  \alpha_3 \tr^2( \Bdis_{,i})  \right) \,.
\end{aligned}
\end{equation}

The optimization results are shown in Figure \ref{fig:results:MM:3} with assuming $\Lc = 1$ mm. The obtained parameters read   $\mu_\textrm{micro} = 5.967$ kN/mm$^2$  , $\mu^*_\textrm{micro} = 392.15$ kN/mm$^2$,  $ \lambda_\textrm{micro} = 10.77$ kN/mm$^2$,  $\mu_c = 808.94$ kN/mm$^2,  \alpha_1  = 0.187 \,, \alpha_2  = 0.318 $ and $  \alpha_3  = 5.65 \,$ . However, the error $r^2$ for the isotropic curvature ($0.000528$ [(kN$\cdot$mm)$^2$]) does not show a significant improvement compared to the simplified curvature ($0.000612$ [(kN$\cdot$mm)$^2$]). Note that the identified parameters ($\Cmicro, \mu_c$) for the two curvature formulations of the micromorphic model exhibit a significant difference, raising doubts about the physical interpretations. The full micromorphic model needs to incorporate the mixed term ($ (\nabla \bu - \Bdis) : \mathbb{C}_\textrm{mixed} : \sym \Bdis$) which will enhance the fitting definitely. 

 \begin{figure}[htpb] 
  \centering
  \fbox{
 \parbox{\textwidth}{ 
We consider forty random modes $i_\textrm{max} = 40$  applied on three sizes $n_\textrm{max}=3$  

$$
{\bar{\bu}_i =   \bB_i \cdot \bx + \bC_i \cdot \bx \otimes \bx  \quad \textrm{on} \quad \partial \B }
$$

where components of $\bB_i$ and $\bC_i$ are randomly generated based on Equations \ref{eq:bcomponents} and \ref{eq:ccomponents}\\

known from previous analysis in Equations \ref{eq:macro} and \ref{eq:mue}: \\

  $\mu_\textrm{macro}  =  5.9 \,  \textrm{ kN/mm$^2$}  \,, \quad 
\mu^*_\textrm{macro}  =   \, 0.627 \,  \textrm{ kN/mm$^2$} \,, \quad \lambda_\textrm{macro}  =  \, 1.748   \textrm{ kN/mm$^2$}\,,  
\quad 
\mu  =   \, 1.537 \,  \textrm{ kN/mm$^2$}  $   \\ 

further assumptions: \\

$\quad \Lc = 1 \, \textrm{mm}, \, \quad  \mathbb{C}_\textrm{mixed} = \bzero$ \\

set initial values:  \\

 $\mu_\textrm{micro} = 26.32 $ kN/mm$^2$ $ \,, \quad  \mu_\textrm{micro}^* =  \, 26.32$ kN/mm$^2$, $\quad \lambda_\textrm{micro} = 51.08 $ kN/mm$^2$,  $\quad \alpha_1 = \alpha_2 = \alpha_3 = 1 $ \\ 

The algorithm delivers:  \\

\begin{tabular}{ccccccccc}
            \toprule
            iteration & $\mu_\textrm{micro}$   & $\mu^*_\textrm{micro}$  &  $\lambda_\textrm{micro}$ &  $\mu_c$ & $\alpha_1$ & $\alpha_2$ & $\alpha_3$ &  $r^2$\\
                &  [kN/mm$^2$] & [kN/mm$^2$]  &  [kN/mm$^2$] &  [kN/mm$^2$]  & & &  & [(kN$\cdot$mm)$^2$]  \\         
            \toprule
               0  & $26.32 $ & $26.32 $  & $ 51.08 $ & $ 1 $& $ 1 $ & $ 1 $ & $ 1 $ & $0.018407$  \\         
            \toprule
               1  & $16.21 $ & $35.89 $  & $ -0.714 $ & $ 1.966 $& $ 0.8$ & $ 0.702 $ & $ 2.254 $ & $0.00180853$  \\         
            \toprule
               2  & $11.76 $ & $62.89 $  & $ 3.425 $ & $ 4.12 $& $ 0.673$ & $ 0.355 $ & $ 3.09 $ & $0.00122141$  \\         
            \toprule
               :  & :  & : & :  & : & : & :  \\   
             \toprule    
               13 & $10.33 $ & $249.7 $  & $ 5.974 $ & $ 189.27 $& $ 0.398$ & $ 0.212 $ & $ 3.69 $ & $0.000956095$  \\         
            \toprule 
               14 & $10.2 $ & $266.56 $  & $ 6.24 $ & $ 523.67 $& $ 0.373$ & $ 0.212 $ & $3.77 $ & $0.000945068$  \\         
            \toprule     
             :  & :  & : & :  & : & : & :  \\   
             \toprule  
               25 & $5.974 $ & $430.96 $  & $ 11.41 $ & $ 773.63 $& $ 0.148$ & $ 0.325 $ & $ 5.883 $ & $0.000534085$  \\         
            \toprule  
              26 & $5.959 $ & $406.16 $  & $ 11.05 $ & $ 784.36 $& $ 0.166$ & $ 0.318 $ & $ 5.758 $ & $0.000532797$  \\               
            \toprule     
             :  & :  & : & :  & : & : & :  \\   
             \toprule  
               49 & $5.966 $ & $392.15 $  & $ 10.77 $ & $ 808.91 $& $ 0.187$ & $ 0.318 $ & $ 5.65 $ & $0.000528028$  \\         
            \toprule  
               50 & $5.967 $ & $392.15 $  & $ 10.77 $ & $ 808.94 $& $ 0.187$ & $ 0.318 $ & $ 5.65 $ & $0.000528028$  \\ 
            \bottomrule
\end{tabular} \\ \\

final parameters set: \\ 

 $\mu_\textrm{micro} = 5.967$ kN/mm$^2$, $\quad \mu^*_\textrm{micro} = 392.15$ kN/mm$^2$,  $\quad \lambda_\textrm{micro} = 10.77$ kN/mm$^2$,   \\
 
  $\mu_c = 808.94$ kN/mm$^2,  \quad \alpha_1  = 0.187 ,  \quad \alpha_2  = 0.318 ,  \quad \alpha_3  = 5.65 \,$

}}
  \caption{Results for the parameter identification algorithm for the full micromorphic model with isotropic curvature. }
    \label{fig:results:MM:3}
\end{figure} 

\subsection{Discussion of results} 
\begin{figure}[htpb]
	\centering
		\begin{subfigure}[b]{0.95\textwidth}
		\includegraphics[width=\textwidth]{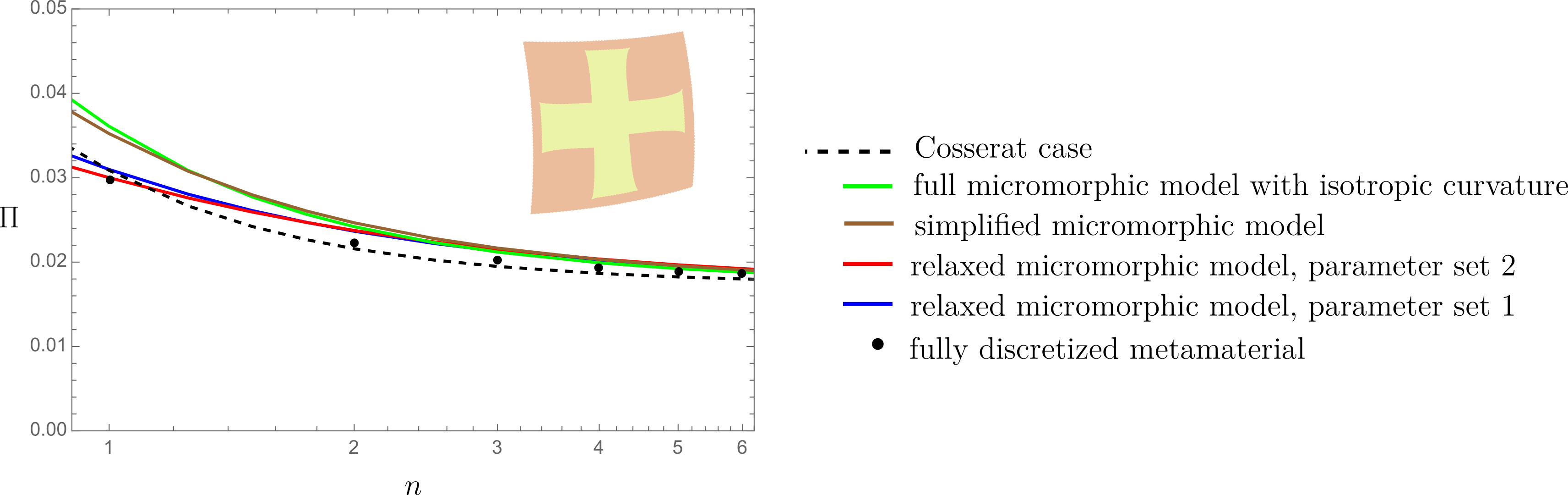}
		\caption*{$i = 1$}
	\end{subfigure}
		\begin{subfigure}[b]{0.49\textwidth}
		\centering
		\includegraphics[width=\textwidth]{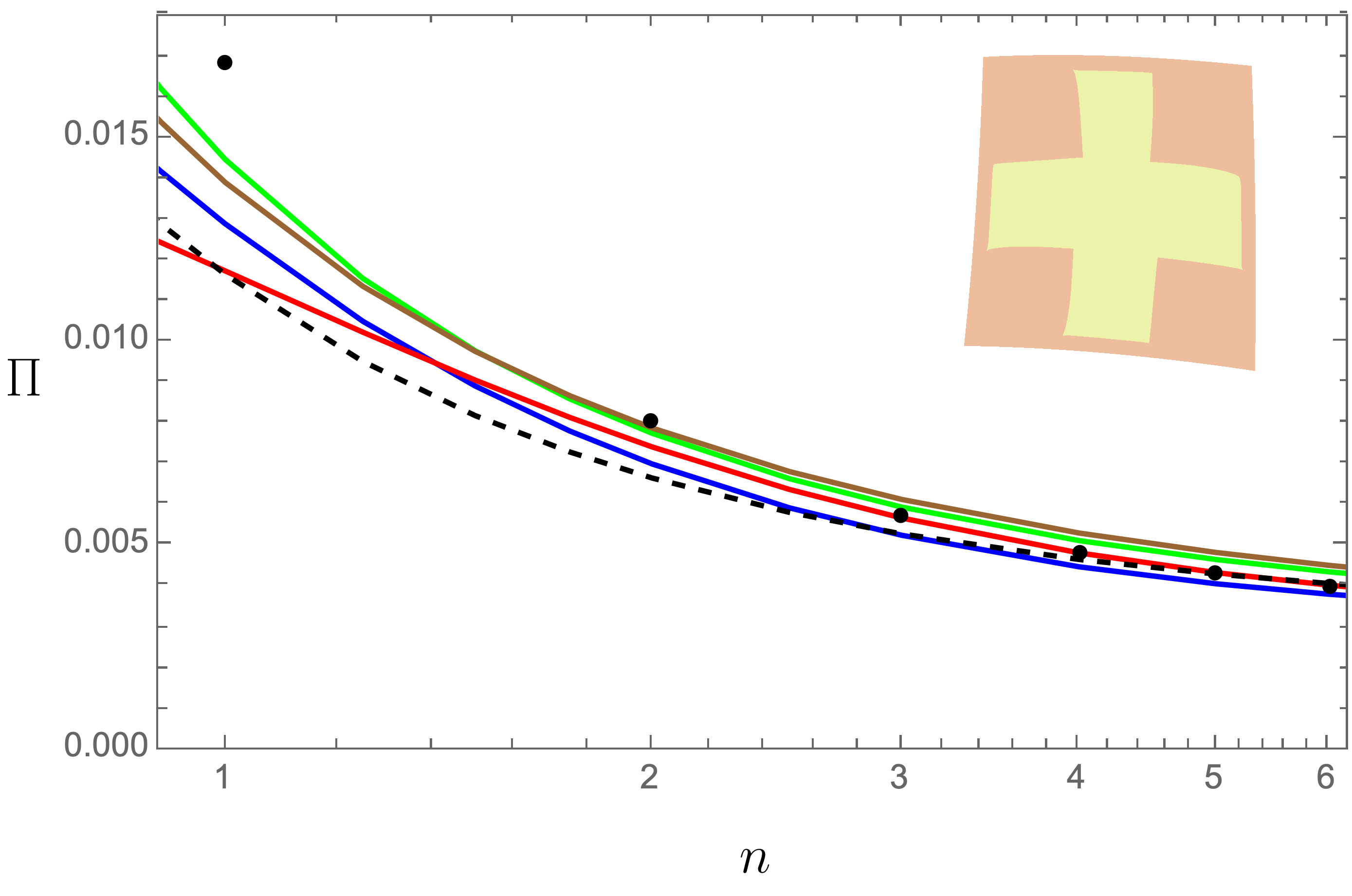}
		\caption*{$i = 5$}
	\end{subfigure}
	\begin{subfigure}[b]{0.49\textwidth}
		\centering
		\includegraphics[width=\textwidth]{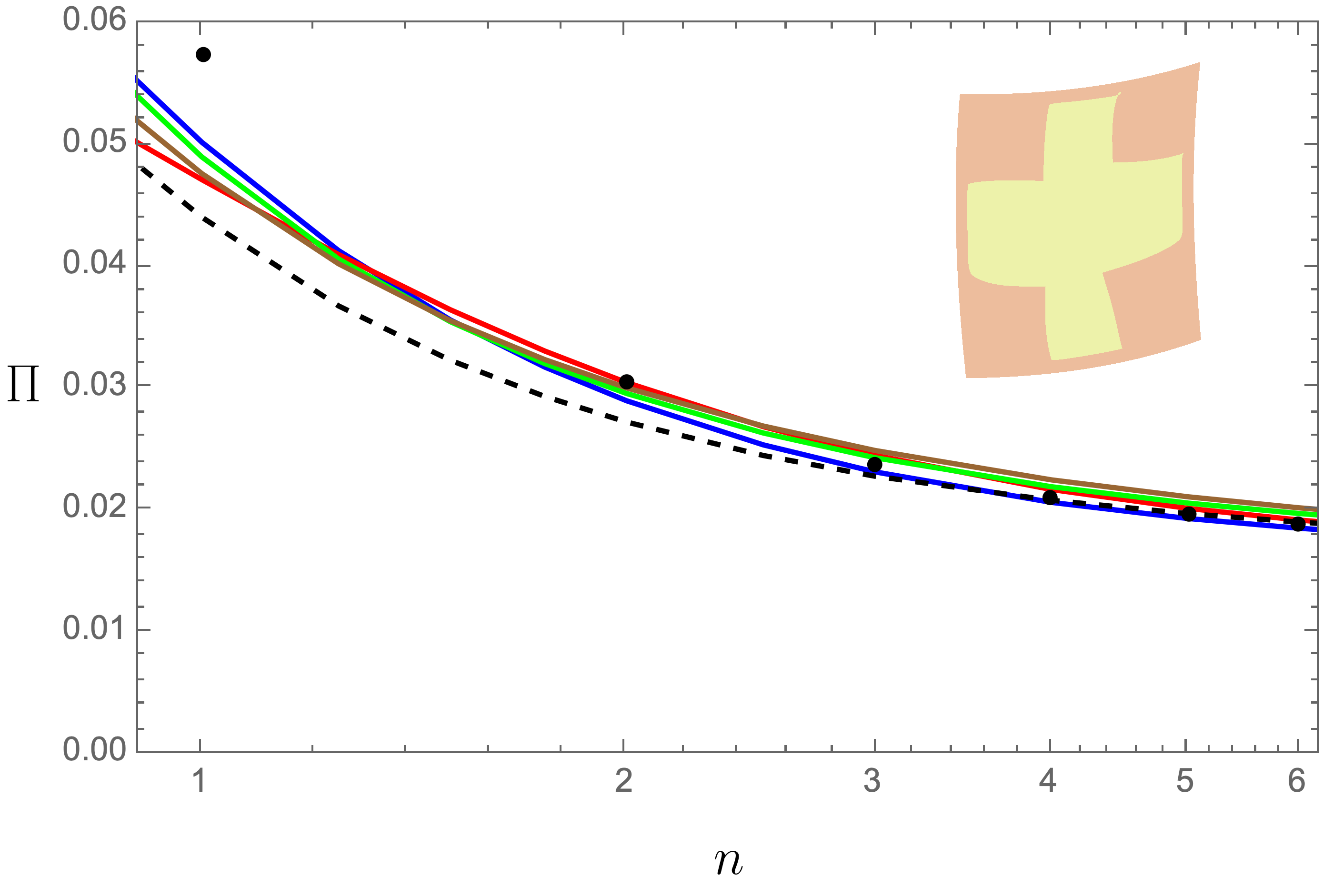}
		\caption*{$i = 10$}
	\end{subfigure}
	\begin{subfigure}[b]{0.49\textwidth}
		\centering
		\includegraphics[width=\textwidth]{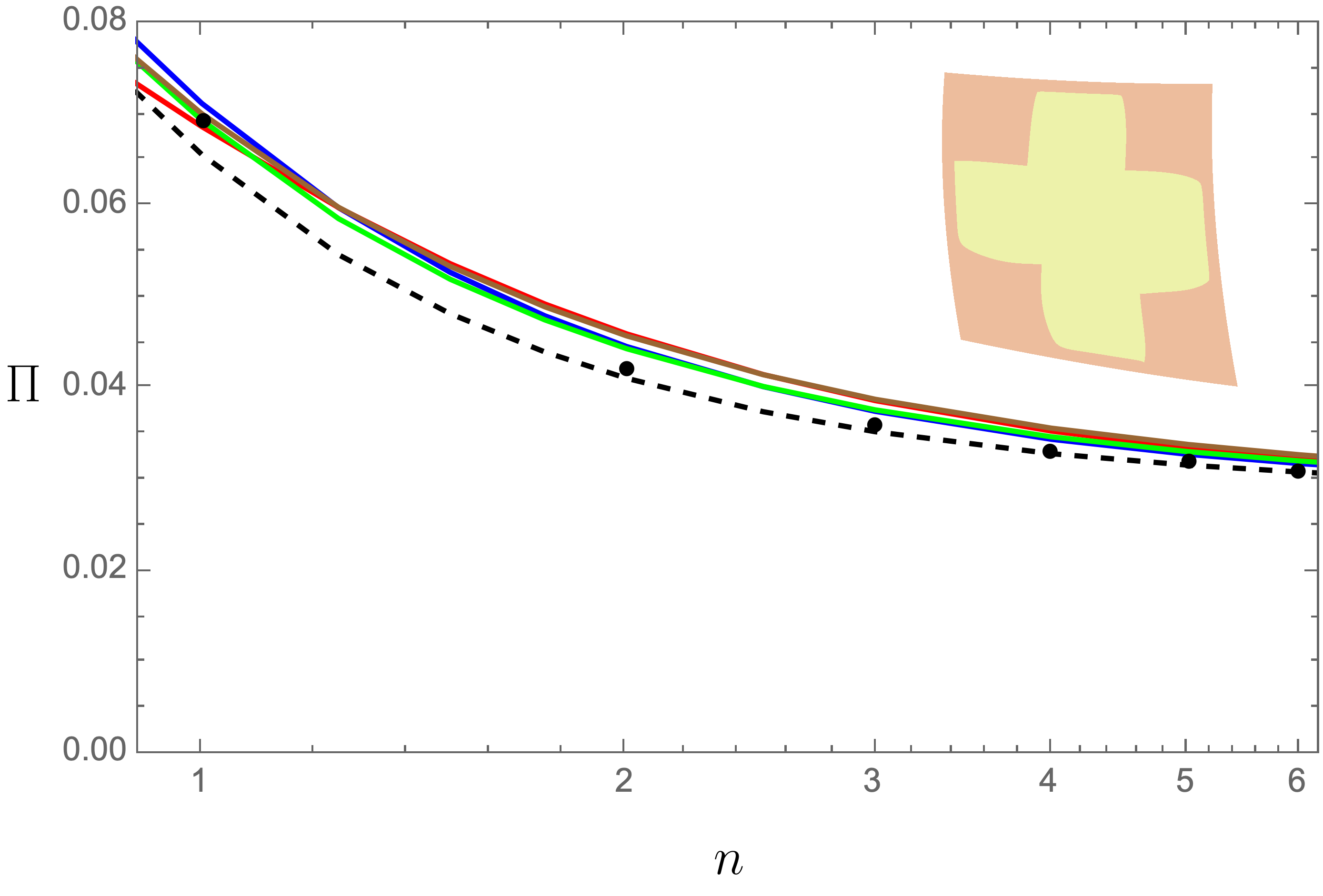}
		\caption*{$i = 15$}
	\end{subfigure}
	\begin{subfigure}[b]{0.49\textwidth}
		\centering
		\includegraphics[width=\textwidth]{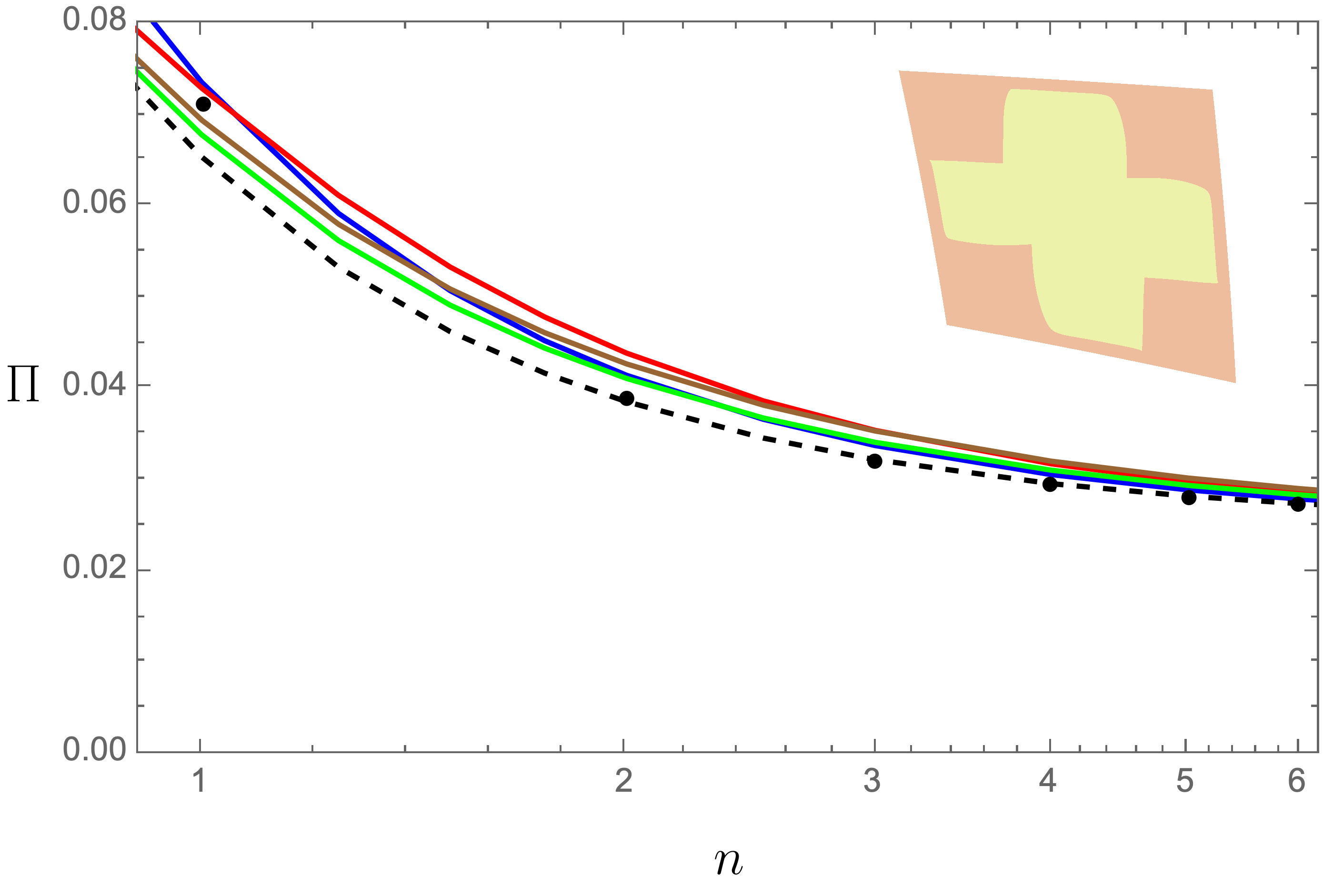}
		\caption*{$i = 20$}
	\end{subfigure}
	\begin{subfigure}[b]{0.49\textwidth}
		\centering
		\includegraphics[width=\textwidth]{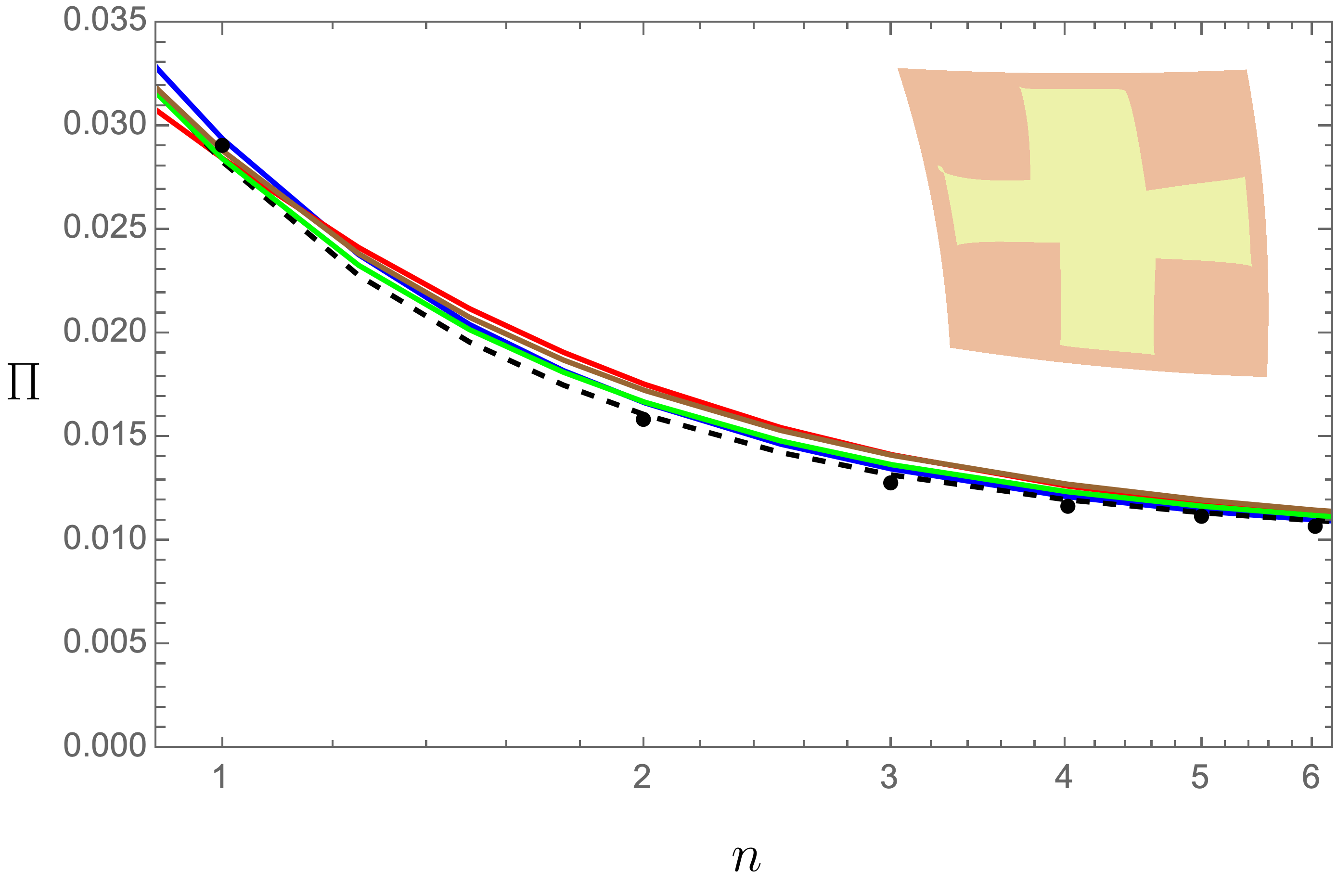}
		\caption*{$i = 30$}
	\end{subfigure}
	\begin{subfigure}[b]{0.49\textwidth}
		\centering
		\includegraphics[width=\textwidth]{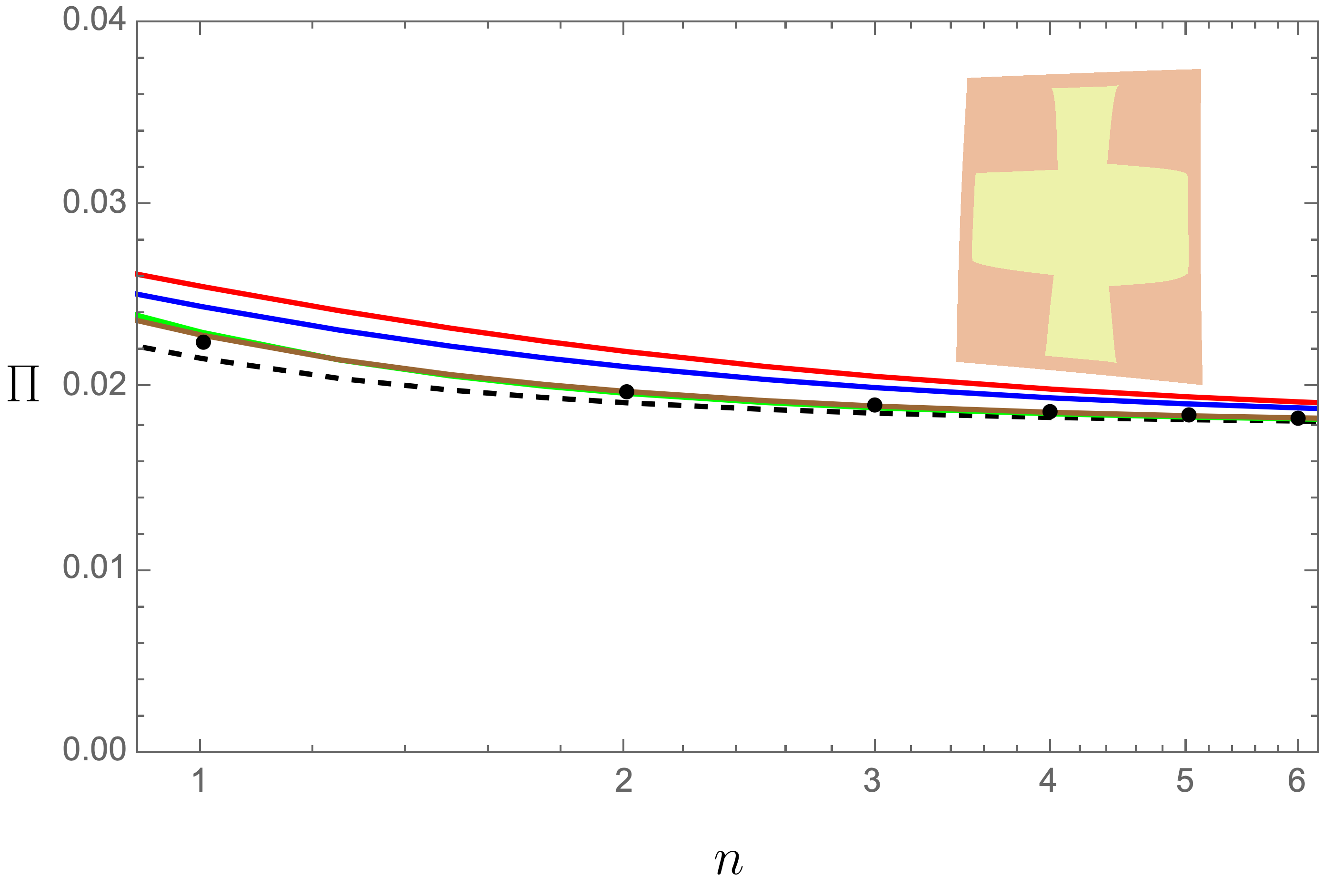}
		\caption*{$i = i_\textrm{max} = 40$}
	\end{subfigure}
	\caption{The total energy of the heterogeneous material, the relaxed micromorphic model,  the Cosserat case, and the  Eringen-Mindlin full micromorphic model with both the simplified and isotropic curvature. We show 7 random deformation modes with the deformed shape for $n=1$ enhanced by a factor of three.   }
		\label{Fig:fitting:6modes}
\end{figure}

The average error of the relaxed micromorphic model concerning 120 reference heterogeneous solutions is $5.3\%$ for parameter set 1 in Figure \ref{fig:results:RMM}, and $7.5\%$ for parameter set 2 in Figure \ref{fig:results:RMM:2}. The Cosserat case leads to an average error $13.7 \%$. The full micromorphic model gives an average error of $7.2 \%$ for the simplified curvature and $6.3 \%$ for the isotropic curvature. Figure \ref{Fig:fitting:6modes} displays the results of fitting 7 deformation modes from the 40 random modes. We extend the results to encompass not just the first 3 sizes ($n=1,2,3$) utilized in the optimization algorithm but also the first 6 sizes ($n=1,..,6$). Beyond these sizes, size-effects are not observed, and standard homogenization theory becomes applicable. Both relaxed micromorphic model and the full micromorphic model demonstrate overall good agreement with the heterogeneous solutions. While the Cosserat model shows the poorest fitting among the examined models, no clear superior model stands out but the unconstrained relaxed micromorphic model (parameter set 1) demonstrates the least error. 

\begin{figure}[htpb]
	\centering

				\includegraphics[width=\textwidth]{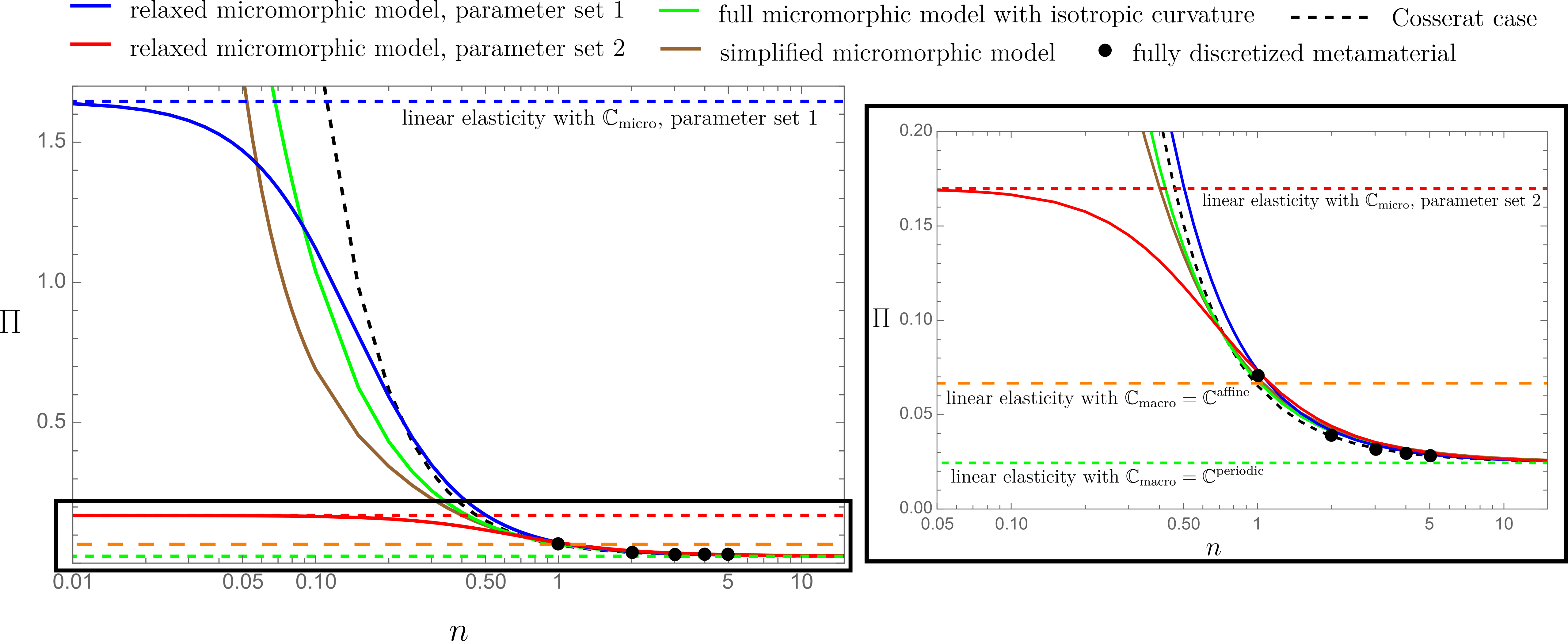}
	\caption{The total energies for the deformation mode $i=20$ varying the size $n$. The relaxed micromorphic model is bounded from above and below by linear elasticity with micro- and macro elasticity tensor, respectively, while the Cosserat case and the simplified full Eringen-Mindlin micromorphic model are not bounded from above. For the heterogeneous case, $n$ must be a natural number, while for the enriched continua $n$ can vary as a real positive number. }
	\label{Fig:fitting:mode_i=20}
\end{figure}

We showcase the characteristics of the relaxed micromorphic model for one deformation mode $i=20$ in Figure 	\ref{Fig:fitting:mode_i=20}. The distinct behavior of the model, functioning as a two-scale linear elasticity model between the macro- and micro-scales, is clearly observed. All the examined models retrieve the results of classical homogenization $\Cmacro$ for large values of $n$ by construction. Parameter set 1 of the relaxed micromorphic model results in a very stiff micro elasticity tensor $\Cmicro$, which is five times stiffer than the homogeneous matrix for the examined deformation mode. In contrast, parameter set 2 yields a softer micro elasticity tensor, nearly half the stiffness of the homogeneous matrix for the same deformation mode.

\begin{figure}[htpb]
	\centering
		\begin{subfigure}[b]{0.49\textwidth}
		\centering
		\includegraphics[width=\textwidth]{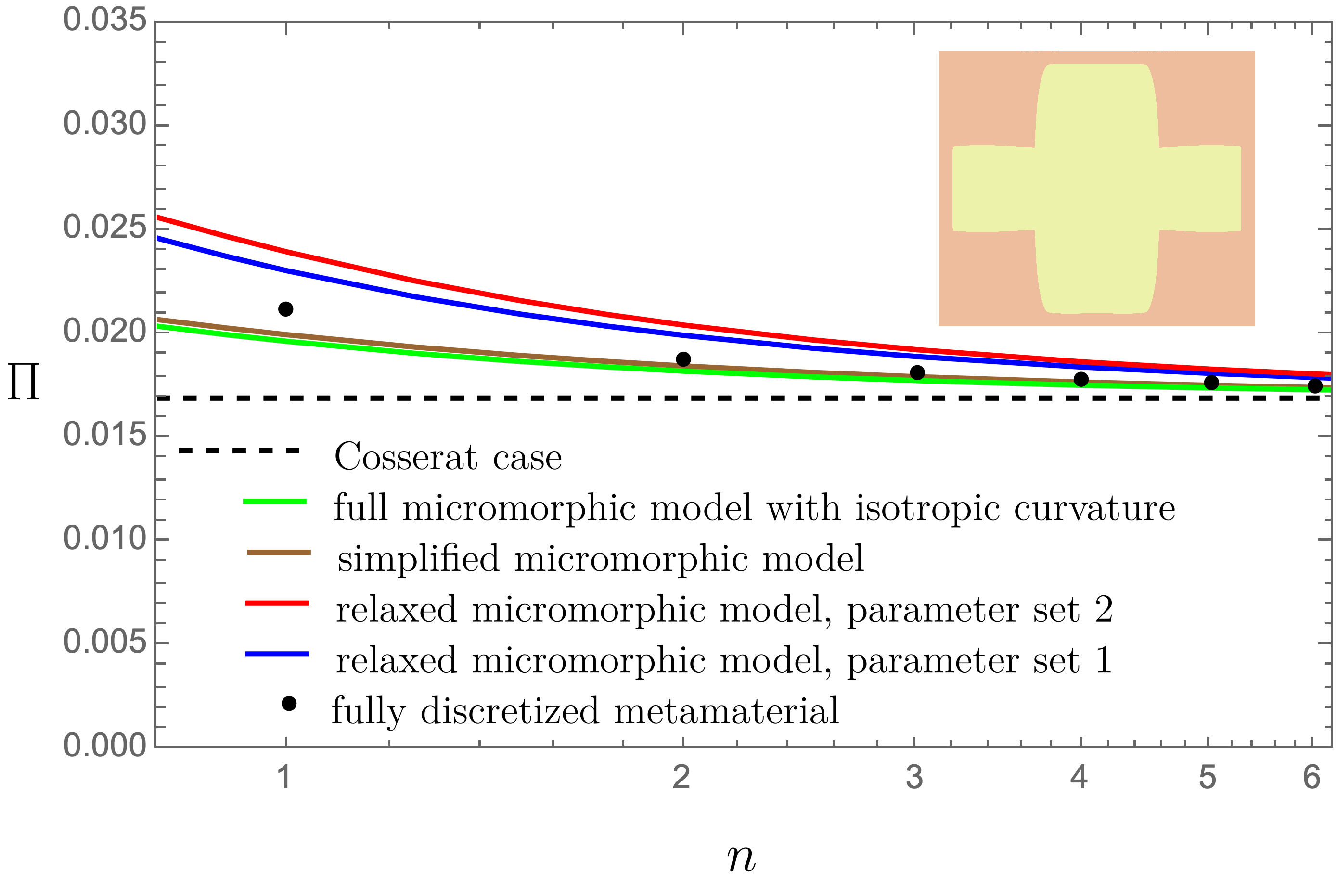}
		\caption*{$\bar\bu^T = (0.05 x,0)$}
	\end{subfigure}
		\begin{subfigure}[b]{0.49\textwidth}
		\centering
		\includegraphics[width=\textwidth]{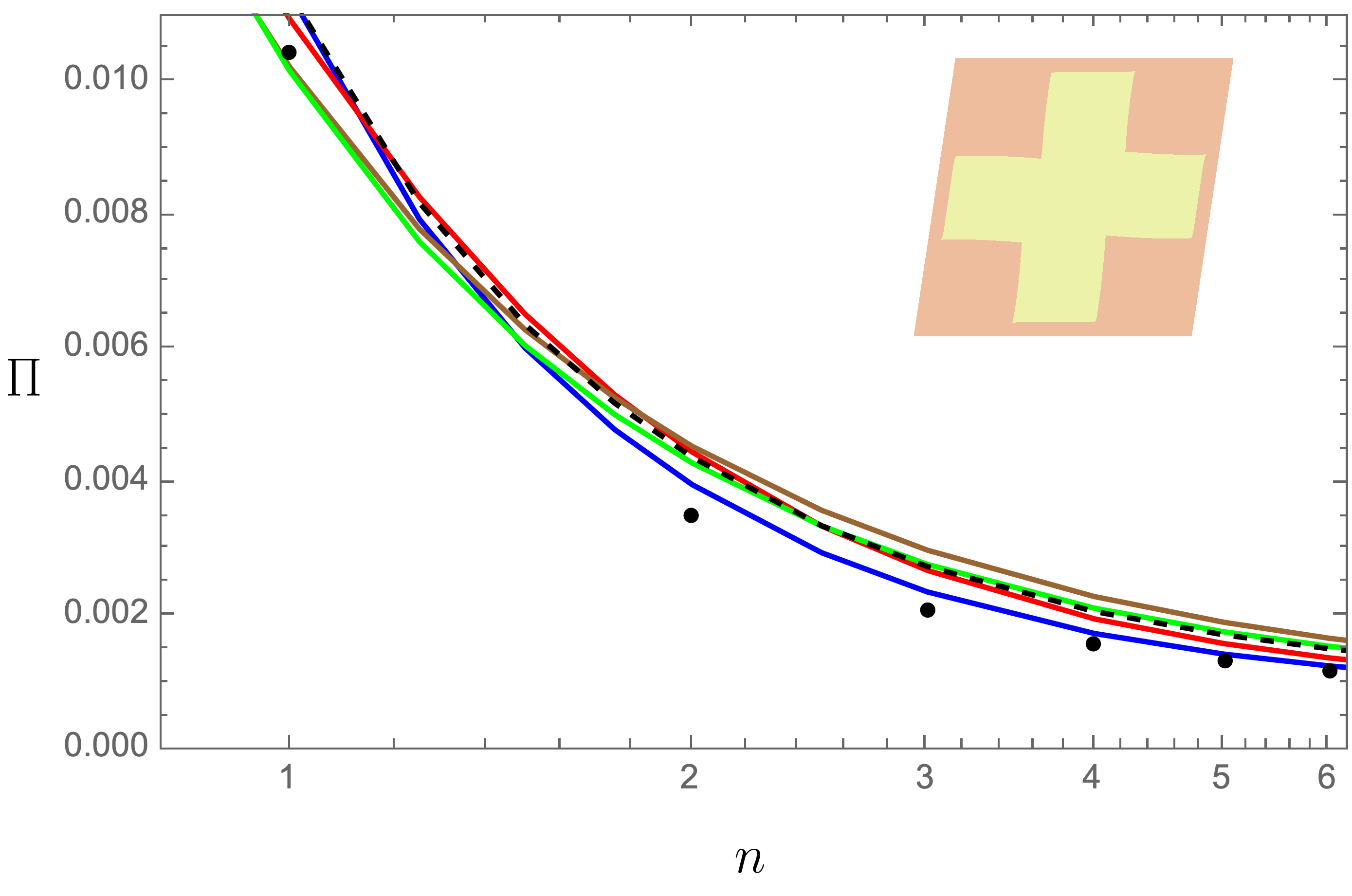}
		\caption*{$\bar\bu^T = (0.05 y,0)$}
	\end{subfigure}
	\begin{subfigure}[b]{0.49\textwidth}
		\centering
		\includegraphics[width=\textwidth]{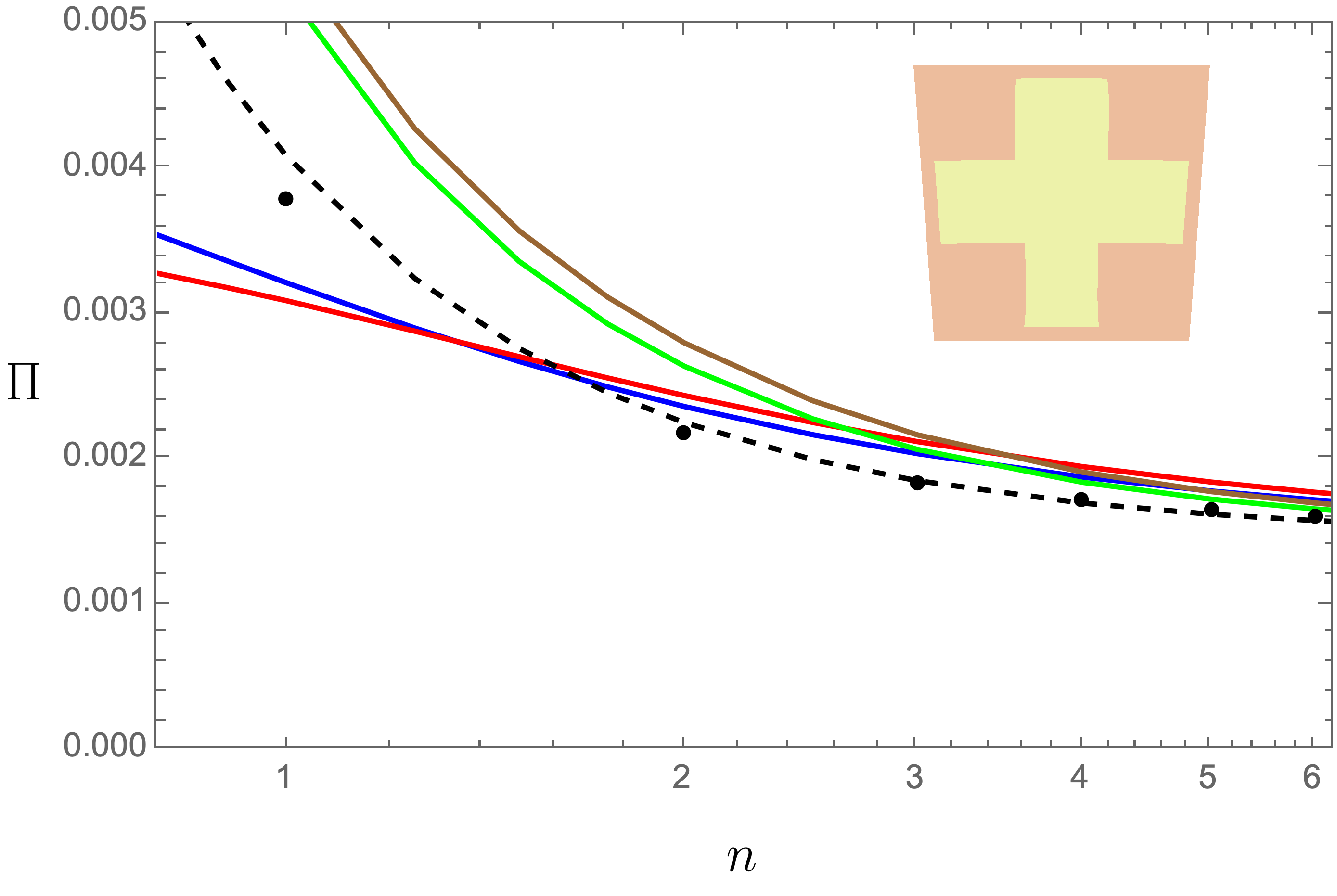}
		\caption*{$\bar\bu^T = (0.05 xy,0)$}
	\end{subfigure}
	\begin{subfigure}[b]{0.49\textwidth}
		\centering
		\includegraphics[width=\textwidth]{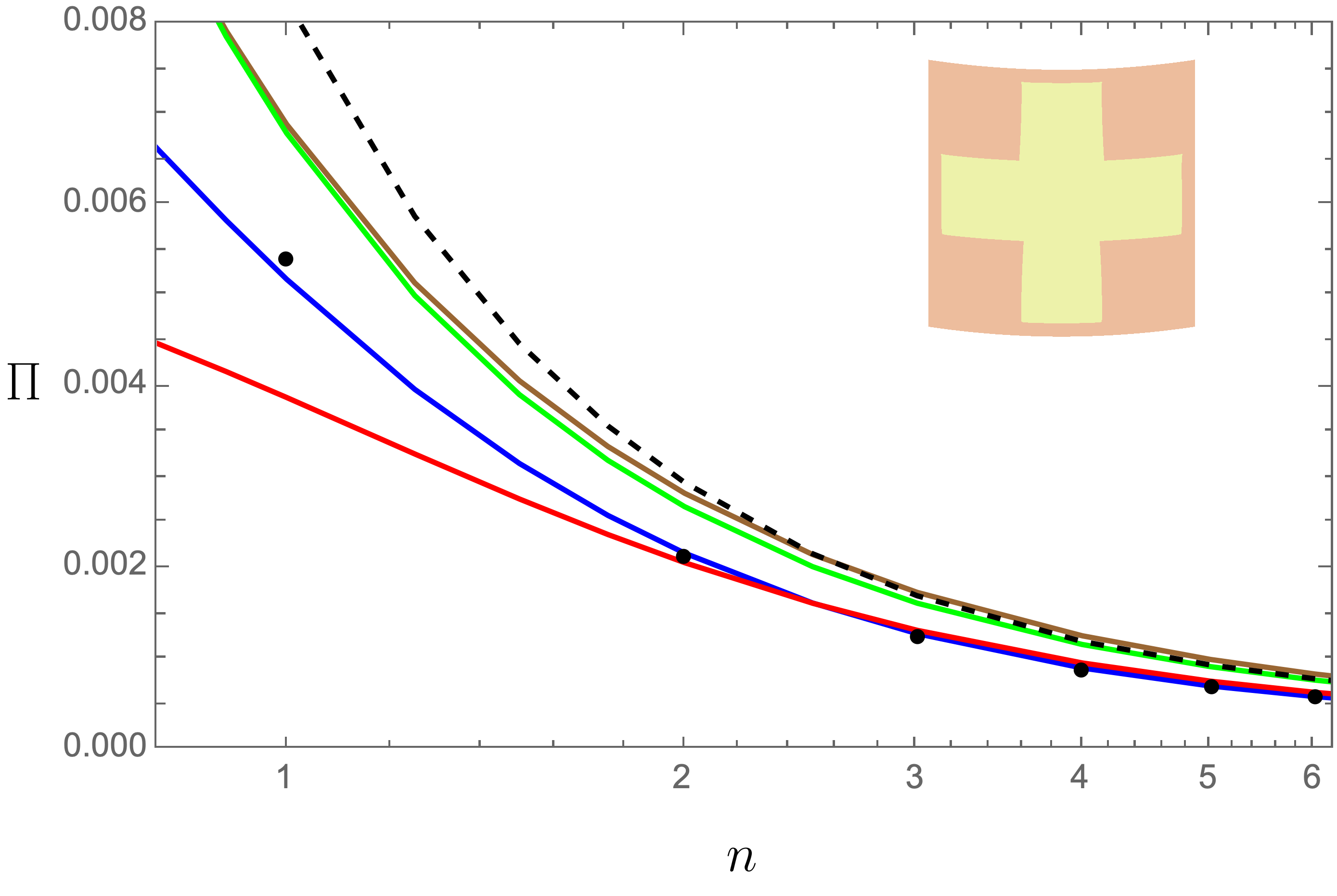}
		\caption*{$\bar\bu^T = (0,0.05 x^2)$}
	\end{subfigure}
	\caption{The total energy of the heterogeneous material and the homogeneous relaxed micromorphic model is examined for the two parameter sets, varying the size for four modes that were not included in the algorithm. We also show the outcomes of the Cosserat and the simplified Eringen-Mindlin micromorphic model. Dirichlet boundary conditions are set on the whole boundary $\bu = \bar\bu \quad \textrm{on} \quad \B$. The deformed shape is shown for $n=1$  enhanced by a factor of three. }
		\label{Fig:fitting:validation}
\end{figure}

We test the obtained parameters of the relaxed micromorphic model in Figure \ref{Fig:fitting:validation} for four different selected deformation modes. Two modes are first-order modes, and the other two modes are second-order modes. The results are satisfactory. Due to the relaxed curvature expression in the RMM, we do not expect to achieve a perfect fit, considering only four unknown parameters to describe the relaxed micromorphic model. The general micromorphic theory already considers 18 parameters in the 3D-isotropic linear elastic case and 14 for the plane strain isotropic case without mixed terms.  The Cosserat model predicts no size-effects for axial loading as expected, and displays good result for shear loading. However, the relaxed micromorphic model still demonstrates better overall agreement. Using the full gradient of the micro-distortion with a single characteristic length parameter does not enhance the overall fitting (the simplified Eringen-Mindlin micromorphic model). We tested whether a better fit can be achieved by using an isotropic curvature with 3 independent parameters. However, this did not lead to a significant improvement. 

Our understanding of the relaxed micromorphic model is reformulated as an extension of classical first-order homogenization to a two-scale elasticity model, not necessarily associated with higher-order deformation modes. It is essential to note that the term $\Cmicro$ appearing in both the relaxed micromorphic model and the classical micromorphic theory does not hold the same interpretation. In the context of the relaxed micromorphic model, it is clearly associated to the bounded stiffness property for small sizes.

\newpage

\section{Conclusions}
\label{Section:Conclusions}

We have successfully established a computational procedure to determine the unknown material parameters of the relaxed micromorphic model. Following a brief consistency check of the presented method for linear elasticity under both affine and periodic boundary condition, which yield correct results in a single iteration, we transferred our methodology to the relaxed micromorphic model. Considering the structural simplicity of the model, only the microscopic elasticity tensor $\Cmicro$ and a scalar associated with curvature remain as unknowns, given that the macroscopic elasticity tensor $\Cmacro$ is known a priori from classical first-order homogenization and  $\Ce$ is uniquely determined once $\Cmicro$ is known.

For our specific cubic unit-cell, only four parameters were included in the algorithm. The algorithm is built on the Hill-Mandel energy equivalence postulate for various deformation modes and sizes, employing a least square fitting approach for implementation. Thus, we have avoided the classical micro-macro transition schemes that rely on the concept of a representative volume element and we eliminated the need to deal with various stress quantities and their microscopic equivalents, which are not known to us yet. Initially, the algorithm produced a good fit, but we found it necessary to introduce an additional criterion ensuring that the microscopic elasticity tensor cannot exceed the stiffness of the homogeneous stiff matrix for any deformation mode. Upon comparing the obtained unknowns with and without this criterion, the results showed overall good agreement.  Additionally, we conducted a comparison between the outcomes of the relaxed micromorphic model and those derived from the Cosserat model, which uses a skew-symmetric micro-distortion field, as well as the simplest case of the full Eringen-Mindlin micromorphic model where we utilized the full gradient of the micro-distortion field as a measurement for the curvature but associated with only a single characteristic length parameter. We attempted to improve the fitting of the full Eringen-Mindlin micromorphic model by employing an isotropic curvature, resulting in a slight enhancement.  The relaxed and full micromorphic continuum exhibit a comparable level of fitting to the fully resolved heterogeneous solution, while the Cosserat continuum performs relatively worse.

We acknowledge that the full Eringen-Mindlin micromorphic model is, of course, capable of exhibiting better fitting than the relaxed micromorphic model. However, for simplified forms of the full micromorphic model with a significantly reduced number of parameters but comparable to the number of parameters of the relaxed micromorphic model or even more, no improvement is achieved. Extending the number of parameters of the full micromorphic model to include the full cubic curvature with ten parameters and considering mixed terms is not feasible for topical engineering applications or optimization procedures. Thus, we can claim that the relaxation of the curvature to considering only the Curl, as in the relaxed micromorphic model, is a reasonable simplification leveraging the strength of the very simple Cosserat model and the more complex full micromorphic model. Comparing the results of the introduced optimization with the outcomes of available micro-macro homogenization schemes will be interesting. However, these micro-macro transition relations for the relaxed micromorphic model are not yet available but are intended to be investigated in future works.

{\bf Acknowledgment} \\
Funded by the Deutsche Forschungsgemeinschaft (DFG, German research Foundation) -  Project number 440935806 (SCHR 570/39, SCHE 2134/1, NE 902/10) within the DFG priority program 2256.

% === list of references
% ------- layout-datei --------------
%\bibliographystyle{abbrv}

\bibliographystyle{plainnat}
%\bibliographystyle{plaindin}
%\bibliographystyle{plaindin_shortname2}
%\bibliographystyle{elsarticle-num-names}
% ------- bib-datei --------------
\FloatBarrier
{\footnotesize
\bibliography{micmag_01}
}
%=== appendix

\end{document}